\documentclass[11pt,leqno,openany]{article}
\usepackage{amsmath,amssymb,amsthm,textcomp,times}
\usepackage[T1]{fontenc}
\usepackage[latin1]{inputenc}
\usepackage{colortbl}
\usepackage[english]{babel}
\usepackage{dsfont}
\newcommand{\bfm}[1]{\boldsymbol{#1}}

\newtheorem{theorem}[subsection]{Theorem}

\newtheorem{definition}[subsection]{Definition}
\newtheorem{lemma}[subsection]{Lemma}
\newtheorem{proposition}[subsection]{Proposition}

\newtheorem{hypothesis}[subsection]{Assumption}

\def \D{ \mathbb D }

\def \H{ \mathbb H}
\def \M{\mathbb M}

\def \R{ \mathbb R }

\def \reel{ \mathbb R }
\def \nat{ \mathds{N}}

\def \ent{\mathds{Z}}

\def \q{ {\rm I}\!\!\!{\rm Q} }
\def \E{ \mathbb  E }
\def \prob{ \mathbb P }
\def \P{ \mathbb P  }

\def \tor{\mathds{T}}
\def \pm{ \bar{\prob}^\varepsilon}
\def \em{ \bar{\E}^\varepsilon}
\begin{document}
\author{Rémi Rhodes\\
Université Paris-Dauphine, Ceremade\\ Place du Maréchal De Lattre De Tassigny
\\75775 Paris cedex 16 - FRANCE,\\e-mail: rhodes@ceremade.dauphine.fr}
\title{Homogenization of locally stationary diffusions with possibly degenerate diffusion matrix}
\maketitle
\begin{quote}
{\bf Abstract :} This paper deals with homogenization of second
order divergence form parabolic operators with locally stationary
coefficients. Roughly speaking, locally stationary coefficients have
two evolution
 scales: both an almost constant microscopic one and a smoothly varying macroscopic one. The homogenization procedure
  aims to give a macroscopic approximation that takes into account the microscopic heterogeneities. This paper follows
 \cite{rhodes:06} and improves this latter work by considering possibly degenerate diffusion matrices.

\vspace{2mm} \noindent {\bf Résumé :} Nous étudions
l'homogénéisation d'opérateurs paraboliques du second ordre sous
forme divergence à coefficients localement stationnaires. Ces
coefficients présentent deux échelles d'évolution: une évolution
microscopique presque constante et une évolution macroscopique
régulière. La théorie de l'homogénéisation consiste à donner une
approximation macroscopique de l'opérateur initial qui tient compte
des hétérogénéités microscopiques. Cet article fait suite à
\cite{rhodes:06} et généralise ce dernier en considérant des
matrices de diffusion pouvant dégénérer.
\end{quote}

%%%%%%%%%%%%%%%%%%%%%%%%%%%%%%%%%%%%%%%%%%%%%%%%%%%%%%%%%%%%%%%%%%%%%%%%%%%%%%%%%%%%%%%%%%%%%%%%%%%%%%%%

\section{Introduction}

%%%%%%%%%%%%%%%%%%%%%%%%%%%%%%%%%%%%%%%%%%%%%%%%%%%%%%%%%%%%%%%%%%%%%%%%%%%%%%%%%%%%%%%%%%%%%%%%%%%%%%%%%%%
This paper follows \cite{rhodes:06} and deals with homogenization of
second order PDEs with locally stationary coefficients by means of
probabilistic tools. More precisely, we aim at describing
 the asymptotic behavior, as $\varepsilon$ goes to
$0$, of the following Stochastic Differential Equation (SDE)
\begin{equation}\label{eq:itroEDS}
  X^\varepsilon_t=x+\frac{1}{\varepsilon}\int_0^tb\big(\omega,\frac{X^\varepsilon_r}{\varepsilon},X^\varepsilon_r\big)\,dr
  +\int_0^tc\big(\omega,\frac{X^\varepsilon_r}{\varepsilon},X^\varepsilon_r\big)\,dr+\int_0^t\sigma\big(\omega,\frac{X^\varepsilon_r}{\varepsilon},X^\varepsilon_r\big)\,dB_r,
\end{equation}
where $B$ is a standard d-dimensional Brownian motion and the
parameter $\omega$ evolves in a random medium $\Omega$, that is a
probability space with suitable stationarity and ergodicity
properties. For each fixed value of the parameter $y\in\R^d$, the
coefficients $b(\omega,\cdot,y) $, $c(\omega,\cdot,y) $ and $
\sigma(\omega,\cdot,y)$ are stationary random fields (the parameter
$\omega$ stands for this randomness). That is why they are said to
be locally stationary. The generator ${\cal L}^\varepsilon$ of the
process $X^\varepsilon$ can be written in divergence form as
\begin{equation}\label{eq:introL}
  {\cal L}^\varepsilon=\frac{1}{2}e^{2V(x)}\sum_{i,j=1}^d\frac{\partial}{\partial
  x_i}\big(e^{-2V(x)}[a+H](\omega,\frac{x}{\varepsilon},x)\frac{\partial}{\partial
  x_j}\,\big)
\end{equation}
for an antisymmetric matrix $H$, a real-valued function $V$ and
$a=\sigma\sigma^*$.

Let us first briefly outline the chronological approach of this
issue. The convergence of the previous SDE (or the connected PDE)
has been first established in the locally periodic case, that is
when the coefficients are deterministic and periodic with respect to
the variable $x/\varepsilon$ \cite{benpard,lions}. Due to the lack
of compactness of a random medium, the random case raises more
difficulties. As far as we know, the first work in this context is
due to Olla and Siri in \cite{siri}. The authors considered a
nearest neighbors random walk on $\ent$ evolving in a locally
stationary environment. They established an invariance principle for
this process under diffusive scaling of space and time. The main
tool of the
 proof is the explicit formula of the correctors, which only holds in dimension one under a strong diffusivity
 condition.

 In \cite{rhodes:06}, an alternative approach is suggested, which is not restricted to the dimension
 one. As in the locally periodic setting, the method is based on a local
 analysis of the microscopic behavior (corresponding to the variable
 $x/\varepsilon$) of the process $X^\varepsilon$ to construct the so-called correctors and to identify the
 limiting process. However, unlike the locally periodic case, these
 correctors turn out to have bad asymptotic properties at a
 macroscopic scale, in the sense that the classical ergodic theory
 cannot describe their asymptotic behavior. Overcoming this issue is the main
 contribution of \cite{rhodes:06}. The main assumption is the uniform ellipticity of the matrix $a$, namely that there exits a constant $M>0$
such that for all $x,y,X\in \reel^d$, $$\frac{1}{M}|X|^2\leq
(a(\omega,x,y)X,X)\leq M |X|^2. $$ This condition is very convenient
for two reasons. From the dynamical angle, it ensures the local
ergodicity of the process $X^\varepsilon$. From the technical angle,
it provides strong estimates of the transition densities of the
process $X^\varepsilon$ as well as regularity properties of its
generator. The control of the process $X^\varepsilon$, in particular
its invariant measure and its tightness, is easily derived from this
assumption.

In this present paper, we intend to improve this latter work by
removing the uniform ellipticity assumption. It is replaced by
 microscopic ergodicity conditions (Assumption \ref{hyp:erg}), which seem not too far from
  being minimal to apply classical ergodic theory and then pass to the limit in
  \eqref{eq:itroEDS}. The class of considered coefficients then
  includes possibly degenerate matrices $a$. In other words, we can
  treat diffusion coefficients $a$ that may reduce to $0$ along some
  directions. Under suitable assumptions, we will prove that the
  process $X^\varepsilon$ converges to the solution $\bar{X} $ of a SDE with
  deterministic coefficients, whose generator can be rewritten in
  divergence form as
\begin{equation}\label{eq:introL}
  \bar{L}=(1/2)e^{2V(x)}\sum_{i,j=1}^d\frac{\partial}{\partial
  x_i}\big(e^{-2V(x)}[\bar{A}+\bar{H}](x)\frac{\partial}{\partial
  x_j}\,\big),
\end{equation}
where the so-called homogenized coefficients $\bar{A}$ and $\bar{H}
$ are respectively symmetric positive and antisymmetric. It is worth
emphasizing that $\overline{A}$ may degenerate, even under strong
non-degeneracy assumptions of the initial diffusion coefficient $a$.
We will prove that the limiting diffusion is trapped in a fixed
subspace of $\reel^d$ and possesses strong diffusivity properties
along this subspace.

 We should finally point out that there are only a few papers
  dealing with possibly degenerate diffusion coefficients in the
  whole literature about probabilistic homogenization of SDEs. In
  the periodic setting, recent advances have been made by Hairer and
  Pardoux in \cite{hairer}. Their approach deeply differs from ours.
  They allow the diffusion to be strongly degenerate in some area of
  the torus, and even to reduce to $0$ over an open domain, provided that the diffusion quickly reaches a strongly
  regularizing area (typically, it satisfies a strong Hörmander type
  condition). Our approach does not allow locally such strong
  degeneracies but does not require any regularizing area. As a
  consequence, we can construct examples that are everywhere
  degenerate. Moreover, the technics used in \cite{hairer} rely on
  the compactness of the torus and cannot be adapted to the random
  setting.

\vspace{2mm} The structure of the paper is the following. In section
2, we introduce all the notations and assumptions. Our results are
stated in Section 4 and an example is given in Section 5. The
construction of the corrector is carried out in Section 6. Section 7
deals with the regularity properties of the process $X^\varepsilon$
such as its invariant measure and the Itô formula. Section 8 is
devoted to establishing the asymptotic properties of the process
$X^\varepsilon$. Section 9 explains the proofs of the homogenization
procedure. The tightness of the process $X^\varepsilon$ is treated
separately in Section 10.
%%%%%%%%%%%%%%%%%%%%%%%%%%%%%%%%%%%%%%%%%%%%%%%%%%%%%%%%%%%%%%%%%%%%%%%%%%%%%%%%%%%%%%%%%%%%%%

\section{Setup and Assumptions}\label{sec:setup}

%%%%%%%%%%%%%%%%%%%%%%%%%%%%%%%%%%%%%%%%%%%%%%%%%%%%%%%%%%%%%%%%%%%%%%%%%%%%%%%%%%%%%%%%%%%%%%

\textit{\textbf{Random medium.}} From now on, $d\geq 1$ is a fixed
integer. Following \cite{jikov}, we introduce the following
\begin{definition}\label{medium}
Let $(\Omega, {\mathcal G},\mu)$ be a probability space and
$\left\{\tau_{x};x\in \reel ^d\right\}$ a group of measure
preserving transformations acting ergodically on $\Omega $:

1) $\forall A\in {\cal G},\forall x\in \reel^d$, $\mu (\tau
_{x}A)=\mu (A)$,

2) If for any $x\in \reel^d$ $\tau _{x}A=A$, then $\mu (A)=0$ or
$1$,

3) For any measurable function ${\bfm g}$ on $(\Omega ,{\cal
G},\mu)$, the function $(x,\omega )\mapsto {\bfm g}(\tau_x \omega)$
is measurable on $(\reel^d\times\Omega ,{\cal B}(\reel^d)\otimes
{\cal G})$.
\end{definition}
The expectation with respect to the random medium is denoted by
$\M$. Denote by $L^2(\Omega)$ the space of square integrable
functions, by $|.|_2$ the corresponding norm and by $(.,.)_2$ the
associated inner product. The operators  defined on $L^2(\Omega )$
by $ T_{x}{\bfm f}(\omega )={\bfm f}(\tau_{x}\omega )$ form a
strongly continuous group of unitary maps in $L^2(\Omega )$.
 For every function ${\bfm f}\in L^2(\Omega )$, let $f(\omega,x )=
{\bfm f}(\tau_{x}\omega )$. Each function ${\bfm f}$ in $ L^2(\Omega
)$ defines in this way a stationary ergodic random field on $\reel
^{d}$. In what follows we will use the bold type to denote an
element ${\bfm f}\in L^2(\Omega )$ and the normal type $f(\omega,x
)$ (or even $f(x)$) to distinguish from the associated stationary
field. The group possesses $d$ generators (throughout this paper,
$e_i$ stands for the i-th vector of the canonical basis of $\R^d$)
\begin{equation} D_i {\bfm g} = \lim_{h \rightarrow 0} \frac{T_{h
e_i} {\bfm g} - {\bfm g}}{h} \ {\rm if \ exists},
\end{equation}
which are closed and densely defined. Setting
\begin{equation}
{\cal C}={\rm Span}\left\{{\bfm g} \star \varphi ;{\bfm g}\in
L^{\infty}(\Omega ),\varphi \in C^\infty _c(\reel ^{d} )\right\}, \
{\rm with} \ {\bfm g} \star \varphi(\omega )=\int_{\reel ^{d}}{\bfm
g}(\tau_{x}\omega )\varphi (x)\,dx,
\end{equation}
the space ${\mathcal C}$ is dense in $L^2(\Omega)$ and ${\cal
C}\subset {\rm Dom}(D_i)$ for all $1 \leq i \leq d$, with $D_i({\bfm
g} \star \varphi)=-{\bfm g} \star
\partial  \varphi/\partial x_i$. If ${\bfm g}\in {\rm Dom}(D_i)$, we also have
$D_i({\bfm g} \star \varphi)= D_i{\bfm g} \star \varphi$. For ${\bfm
f}\in \bigcap_{i=1}^d{\rm Dom}(D_i)$, we define the divergence
operator ${\rm Div}$ by ${\rm Div}{\bfm f}=\sum_{i=1}^dD_i{\bfm f}$.
We distinguish this latter operator from the usual divergence
operator on $\reel^d$ denoted by the small type $\textrm{div}$.

\vspace{2mm} \textit{\textbf{Locally stationary random fields.}}
Following the notations introduced just above, for
 a measurable function ${\bfm f}:\Omega\times\R^d\rightarrow \R^n$, ($n\geq 1$),
 we can consider the associated locally stationary random field
 $(x,y)\mapsto {\bf f}(\tau_x\omega,y)=f(\omega,x,y)$ (or even $f(x,y)$).

\vspace{3mm} \textit{\textbf{Structure of the coefficients}}. The
coefficients ${\bfm \sigma}:\Omega\times \reel^d\rightarrow
\reel^{d\times d}, {\bfm H}:\Omega\times \reel^d\rightarrow
\reel^{d\times d}$, $\tilde{{\bfm \sigma}}:\Omega\rightarrow
\reel^{d\times d}$ and $V:\reel^d\rightarrow \reel$ denote
measurable functions with respect to the underlying product
$\sigma$-fields. As explained above, ${\bfm \sigma}$ and ${\bfm H}$
define locally stationary random fields and $\tilde{{\bfm \sigma}}$
a stationary random field. ${\bfm H}$ is antisymmetric. We define
two new matrix-valued functions by ${\bfm a}={\bfm \sigma}{\bfm
\sigma}^*$ and $\tilde{{\bfm a}}=\tilde{{\bfm \sigma}}\tilde{{\bfm
\sigma}}^*$. Furthermore, for some positive constant $\Lambda$, the
coefficients ${\bfm \sigma}$, ${\bfm H}$, $\tilde{{\bfm \sigma}}$
and $V$ satisfy
\begin{hypothesis}{{\bf (Regularity).}}\label{hyp:reg}
For each fixed $\omega\in\Omega$, the coefficients
$\sigma(\omega,.,.)$, $H(\omega,.,.)$ and $\tilde{\sigma}(\omega,.)$
are two times continuously differentiable with respect to each
variable and are, as well as their derivatives up to order two,
$\Lambda$-Lipschitzian and bounded by $\Lambda$. $V$ is three times
continuously differentiable and is, as well as its derivatives up to
order three, bounded by $\Lambda$ and $\Lambda$-Lipschitzian.
\end{hypothesis}

Let us now describe the degeneracies of the matrix ${\bfm a} $.
Roughly speaking, the degeneracies of ${\bfm a}$ are assumed to be
controlled by the reference matrix $\tilde{{\bfm a}}$. To be more
explicit, let us first introduce the
\begin{definition}
Given a $d\times d$ matrix-valued function $g:\R^d\rightarrow
\R^{d\times d}$, a $d\times d $ symmetric matrix $A$ and a real
$C>0$, $g$ is said to be $(C,A)$-controlled if $\forall y,y'\in\R^d$
$$|g(y)|\leq C A, \quad  \text{ and } |g(y)-g(y')|\leq C A|y-y'|,$$
where $|M|=(MM^*)^{1/2}$ stands for the absolute value of the matrix
$M$ (given 2 symmetric matrices $A,B$, the relation $A\leq B$ means
that the matrix $B-A$ is symmetric positive).
\end{definition}
We now precise the control of ${\bfm a}$ by $\tilde{{\bfm a}}$:
\begin{hypothesis}\label{hyp:cont}{{\bf (Control).}} We
assume that $$M^{-1}\widetilde{{\bfm a}}(\omega)\leq {\bfm
a}(\omega,y)\leq M \widetilde{{\bfm a}}(\omega)$$ for some strictly
positive constant $M$ and for every $(\omega,y)\in\Omega\times
\reel^d$. Moreover, for any $i,j\in \{1,\dots,d\}$ and
$(\omega,y)\in\Omega\times \reel^d$, the matrices $\partial
_{y_i}{\bfm a}(\omega,y)$, $\partial^2 _{y_iy_j}{\bfm a}(\omega,y)$,
${\bfm H}(\omega,y) $, $\partial _{y_i}{\bfm H}(\omega,y)$,
$\partial^2 _{y_iy_j}{\bfm H}(\omega,y)$ are $(M,\tilde{{\bfm
a}}(\omega))$-controlled. We further assume that
$$|{\bfm \sigma }(\omega ,y+h)-{\bfm \sigma }(\omega ,y))|^2 \leq
M\widetilde{{\bfm a}}(\omega)|h|^2$$ for any $y,h\in \reel^d$ and
that $\int_{\R^d}e^{-2V(y)}\,dy=1$.
\end{hypothesis}
To ensure the local ergodicity of the process $X^\varepsilon$, we
make the following assumption:
\begin{hypothesis}[{\bf Ergodicity}]\label{hyp:erg}
Let us consider the Friedrich extension (see \cite[p. 53]{fuku} or
Section 5) of the symmetric operator $\tilde{\bfm S}$ defined on
${\mathcal C}\subset L^2(\Omega)$ by $\tilde{\bfm S}= (1/2) \sum_{i,j=1}^d$ $
D_i(\tilde{\bfm a}_{i,j} D_j )$. This extension, still denoted
$\tilde{\bfm S}$, is self-adjoint. We then assume that the
semi-group generated by $\tilde{\bfm S}$ is ergodic, that is its
invariant functions are $\mu$ almost surely constant (see e.g.
Rhodes \cite{rhodes:04}).
\end{hypothesis}

\noindent {\small \textit{\textbf{Remark.} Assumptions \ref{hyp:reg}
may appear restrictive and can surely be relaxed (see \cite{delarue}
for results in this direction in the context of quasilinear PDEs).
In particular, the statement of the homogenization property only
involves the derivatives of order one with respect to $y\in\R^d$
(see Theorem \ref{maintheorem}). However, it avoids dealing with
heavy regularizing procedures that are not the purpose of this
work.}}

\vspace{2mm}

\textit{\textbf{Diffusion in a locally ergodic environment.}} For
$j=1,\dots,d$, we define the coefficients
\begin{equation}\label{coefb}
%\begin{align}
{\bfm  b}_j(\omega ,y)=  \frac{1}{2}\sum_{i=1}^dD_i({\bfm a}+{\bfm
H})_{ij}(\omega,y),\quad {\bfm c}_j(\omega,y)=  \frac{e^{2V(y)}}{2}
 \sum_{i=1}^d \partial _{y_i}\big(e^{-2V}[{\bfm a}+{\bfm
 H}]_{ij}\big)(\omega,y).
 \end{equation}
From Assumption \ref{hyp:reg}, the functions $b_j(\omega ,.,.) $ and
$c_j(\omega ,.,.) $ are Lipschitzian so that, for a starting point
$x\in\reel^d$ and $\varepsilon>0$, we can consider the strong
solution $X^{\varepsilon}$ of the following Stochastic Differential
Equation (SDE) with locally stationary coefficients:
\begin{equation}\label{eq:diff} X^{\varepsilon
}_t=x+\frac{1}{\varepsilon}\int_{0}^{t}b
\left(\overline{X}^{\varepsilon }_r, X^{\varepsilon }_r \right)\,dr+
\int_{0}^{t}c \left(\overline{X}^{\varepsilon }_r, X^{\varepsilon
}_r \right)\,dr+\int_{0}^{t}\sigma \left(
\overline{X}^{\varepsilon}_r,X^{\varepsilon }_r\right)\,dB_r,
\end{equation}
where we have set $\overline{X}^{\varepsilon }_t\equiv
X^{\varepsilon }_t/\varepsilon $ and $B$ is a standard d-dimensional
Brownian motion (the random medium and the Brownian motion are
independent). We point out that the generator of this diffusion
could be formally written in divergence form as
\begin{equation}\label{generator}
{\cal
L}^{\varepsilon}=\frac{1}{2}e^{2V(x)}\sum_{i,j=1}^d\frac{\partial}{\partial
x_i}\big(e^{-2V(x)}[a+H](\omega,x/\varepsilon,x)\frac{\partial}{\partial
x_j}\,\big).
\end{equation}

\noindent \textit{{\small \textbf{Notations.} For the sake of
simplicity, we indicate the starting point $x$ of $X^{\varepsilon}$
by writing, when necessary, $\P^\varepsilon_x$ (and
$\E^\varepsilon_x$ for the corresponding expectation), this avoids
heavy notations as $X^{\varepsilon,x}$. We can then consider the
probability measure $\pm\equiv
\M\int_{\reel^d}\prob^\varepsilon_x[.]e^{-2V(x)}\,dx$ and its
expectation $\em$.
%In the sequel, the generic notations ``$C$'' and
%``$D$'' stand for constants that only refer to $M$ and $\Lambda$.
%Dependencies on additional parameters are always mentioned.
}}
%%%%%%%%%%%%%%%%%%%%%%%%%%%%%%%%%%%%%%%%%%%%%%%%%%%%%%%%%%%%%%%%%%%%%%%%%%%%%%%%%%%%%%%%%%%%%%%%%%%%%%%%
\section{Main Results}
%%%%%%%%%%%%%%%%%%%%%%%%%%%%%%%%%%%%%%%%%%%%%%%%%%%%%%%%%%%%%%%%%%%%%%%%%%%%%%%%%%%%%%%%%%%%%%%%%%%%%%%%%%%
Let us now state the main result of this paper. Under the previous
assumptions, we can prove

\begin{theorem}\label{maintheorem}{{\bf Homogenization.}}
The law $\pm$ of the process
$X^\varepsilon$ weakly converges in $C([0,T];\reel^d)$ towards the
law of the process $X$ that solves the following SDE with
deterministic coefficients (they do not depend on the medium
$\Omega$):
\begin{equation}\label{martprob}
  X_t=
  x+\int_0^t\overline{B}(X_r)\,dr+\int_0^t\overline{A}^{1/2}(X_r)\,dB_r.
\end{equation}
The coefficients $\overline{A}$ and $\overline{B}$ are of class
$C^2$ and are defined, for $y\in \R^d$, by
\begin{subequations}
\begin{align}
\overline{A}(y) & =\lim_{\lambda\to 0}\,\M[(I+{\bfm D}{\bfm
u}_\lambda)^*{\bfm a}(I+{\bfm D}{\bfm
u}_\lambda)(.,y)],\label{hom:a}
\\\overline{H}(y) & =\lim_{\lambda\to 0}\,\M[(I+{\bfm D}{\bfm
u}_\lambda)^*{\bfm H}(I+{\bfm D}{\bfm u}_\lambda)(.,y)],\label{hom:b}\\
\overline{B}(y) & =
(1/2)e^{2V(y)}\partial_y(e^{-2V}[\overline{A}+\overline{H}])(y).\label{hom:c}
\end{align}
\end{subequations}
Formally speaking, for each $y\in\R^d$ and $\lambda>0$, the entries
$\big({\bfm u}^i_\lambda(.,y)\big)_{1\leq i \leq d}$ of the function
${\bfm u}_\lambda(.,y):\Omega\rightarrow\R^d$ solve the following
so-called {\it auxiliary problems}, which are stated on the random
medium
$$ \lambda {\bfm u}^i_\lambda(.,y)-\frac{1}{2}\sum_{j,k}D_j\big[({\bfm a}_{jk}+{\bfm H}_{jk})D_k{\bfm u}^i_\lambda(.,y)\big]={\bfm b}_i(.,y).$$
\end{theorem}

\noindent \textit{{\small \textbf{Remark}. A rigorous description of
${\bfm u}_\lambda(.,y)$ is given in Section \ref{sec:auxiliary}. In
particular, in this degenerate framework, the "gradients" ${\bfm
D}{\bfm u}_\lambda$ do not exist but along the direction
$\tilde{{\bfm \sigma}}$, that is the only expression $\tilde{{\bfm
\sigma}}^*{\bfm D}{\bfm u}_\lambda $ can be given a rigorous sense.
Because of the control of ${\bfm a}$ and ${\bfm H}$ by $\tilde{{\bfm
a}}$ (Assumption \ref{hyp:cont}), it then makes sense to consider
formulae \eqref{hom:a} and \eqref{hom:b} (see Section
\ref{sec:auxiliary} for further details).}}

\vspace{1mm} Since the diffusion coefficient ${\bfm a}$ is allowed
to degenerate, the reader may wonder whether the homogenized
diffusion coefficient may also degenerate. The following proposition
details the structure of the limiting diffusion coefficient
$\bar{A}$:
\begin{proposition}\label{propgeom}{{\bf Geometry of the homogenized coefficients.}}
 The kernel $K={\rm Ker}(\bar{A}(y))$ of $\bar{A}(y)$ does not depend on the point $y\in\R^d$ where it is computed.
For each $y\in\R^d$, $ \overline{B}(y)\in K^{\bot}$ ($K^{\bot}$ is
the orthogonal complement to $K$) and there exists a
  constant $\alpha_{\ref{propgeom}}>0$, such that
\begin{equation*}
\forall y \in \R^d, \ \forall x \in K^{\bot}, \
\alpha_{\ref{propgeom}}^{-1} |x|^2 \langle x, \bar{A}(y) x \rangle
\leq \alpha_{\ref{propgeom}} |x|^2.
\end{equation*}
In other words, for each starting point $x\in\R^d$, the limiting
process $X$ (see \eqref{martprob}) can be seen as the solution of a
SDE defined on $x+K^{\bot}$ with a uniformly elliptic diffusion
matrix $\bar{A}$.
\end{proposition}
%%%%%%%%%%%%%%%%%%%%%%%%%%%%%%%%%%%%%%%%%%%%%%%%%%%%%%%%%%%%%%%%%%%%%%%%%%%%%%%%%%%%%%%%%%%%%%%%%%%%%%%%

\section{Example}\label{example}

%%%%%%%%%%%%%%%%%%%%%%%%%%%%%%%%%%%%%%%%%%%%%%%%%%%%%%%%%%%%%%%%%%%%%%%%%%%%%%%%%%%%%%%%%%%%%%%%%%%%%%%%%%%

Let us consider a simple example in the two dimensional
$2\pi$-periodic case. The 2-dimensional torus $\tor^2$ is seen as the random medium equipped with the induced Lebesgue measure, still denoted by $\mu$ to stick with the notations of the paper . We aim at constructing a degenerate
homogenized coefficient. For this purpose, let us first define
$$\forall x\in\reel^2,\quad \widetilde{\bfm{\sigma}} (x)=\left(\begin{array}{cc}
 1 &  1/c\\
c      &  1 \end{array}\right),$$ where $c\not\in \pi\q$ is a
constant, and $\widetilde{\bfm{a}}=\widetilde{\bfm{\sigma
}}\widetilde{{\bfm \sigma}}^*$. Choose now any smooth function
${\bfm U}:\reel^2\times\reel^2\rightarrow \reel^{2\times 2}$, with
bounded derivatives up to order $2$, $2\pi$-periodic with respect to
its first argument $x\in\R^2$ and satisfying
$$\forall(x,y)\in\reel^2\times \reel^2,\quad
M^{-1}\mathrm{Id}\leq {\bfm U}{\bfm U}^* (x,y)\leq M\,
\mathrm{Id}.$$ Define $\forall (x,y)\in \reel^2\times\reel^2$,
$V(y)= e^{-|y|^2}/\pi$, ${\bfm \sigma} (x,y)=\widetilde{{\bfm
\sigma}}{\bfm U}(x,y)$ and ${\bfm H}=0$. Let us check that these
coefficients satisfy all our assumptions. From the smoothness of the
coefficients, it is plain to see that Assumptions \ref{hyp:cont} and
\ref{hyp:reg} are fulfilled. Assumption \ref{hyp:erg} results from
the Weyl equipartition theorem ($c\not\in\pi\q$). Theorem
\ref{maintheorem} thus holds.

Let us now prove that $\bar{A}$ is degenerate and does not trivially reduce to $0$.  Let us denote by $\widetilde{A}$ the homogenized
coefficient associated to $\tilde{{\bfm a}}$.  From the proof of
Proposition \ref{propgeom}, for any $y\in\R^2$ and $X\in \R^2$, we have
$$C^{-1}\big<X,\widetilde{A}X\big>\leq\big<X,\bar{A}(y)X\big>\leq C
\big<X,\widetilde{A}X\big>=0.$$
So we just have to compute $\widetilde{A}$. Since $\widetilde{\bfm{\sigma}}$ is constant, it is straightforward to check that $\widetilde{A}$ actually matches $\widetilde{\bfm{\sigma}}\widetilde{\bfm{\sigma}}^* $ with the help of \eqref{varform}.  Indeed, for a given smooth function $ {\bfm \varphi}$ defined on $\tor^2$ and $x\in \R^2$, the right-hand side of \eqref{varform} expands as 
\begin{align*}
\M[|\widetilde{\bfm{\sigma}}(D{\bfm \varphi}+x)|^2]=& \M[|\widetilde{\bfm{\sigma}}^*D{\bfm \varphi}|^2]+2\left<\widetilde{\bfm{\sigma}}^*x,\widetilde{\bfm{\sigma}}^*\M[D{\bfm \varphi}]\right>+\left<\widetilde{\bfm{\sigma}}^*x,\widetilde{\bfm{\sigma}}^*x\right>\\ =&\M[|\widetilde{\bfm{\sigma}}^*D{\bfm \varphi}|^2]+\left<\widetilde{\bfm{\sigma}}^*x,\widetilde{\bfm{\sigma}}^*x\right>  .
\end{align*}
The infimum is then clearly reached for ${\bfm \varphi}=0 $.

Finally, we let the reader check that $\widetilde{A}=\widetilde{\bfm{\sigma}}\widetilde{\bfm{\sigma}}^*$ does not reduce to $0$ and that the vector $X_K=[1 \,\,\,-c]^*$ satisfies
$\widetilde{A}X_K=0$. 

In a general way, because of the various geometries of random media, it is not clear whether $\bar{A}$ is degenerate or not. The reader may find in \cite{delarue} examples (in a slightly different framework) where the diffusion matrix reduces to $0$ though the diffusion coefficient ${\bfm \sigma}$ is elliptic over a set of full Lebesgue measure, and conversely, an example where ${\bfm \sigma}$ degenerates and $\bar{A}$ is uniformly elliptic.
\qed
%%%%%%%%%%%%%%%%%%%%%%%%%%%%%%%%%%%%%%%%%%%%%%%%%%%%%%%%%%%%%%%%%%%%%%%%%%%%%%%%%%%%%%%%%%%%%%%%%%%%%%%

\section{Construction of unbounded operators}\label{sec:construction}

%%%%%%%%%%%%%%%%%%%%%%%%%%%%%%%%%%%%%%%%%%%%%%%%%%%%%%%%%%%%%%%%%%%%%%%%%%%%%%%%%%%%%%%%%%%%%%%%%%%%%%%%%%%
Throughout this paper, we will need to construct suitable extensions
of unbounded operators defined on a dense subspace of a given
$L^2$-space. This construction is always the same and follows
\cite[Ch. 3, Sect 3.]{fuku} or \cite[Ch. 1, Sect 2.]{ma}, to which
the reader is referred for further details than those given below.
That is the reason why we explain it in a generic way. We also point
out that the Friedrich extension of $\tilde{{\bfm S}}$ (see
Assumption \ref{hyp:erg}) corresponds to this construction.

Consider a probability space $\Omega$ equipped with a probability
measure $\P$, a dense subspace ${\cal D}$ of $L^2(\Omega;\P)$, a
 positive symmetric bilinear form $\left<\cdot,\cdot\right>$ defined
on ${\cal D}\times {\cal D}$ ($\|\cdot\|$ denotes the corresponding
semi-norm) and a bilinear form $B$ on ${\cal D}\times {\cal D}$ that
satisfies for any $\varphi,\psi\in {\cal D} $
\begin{equation}\label{equivB}
\alpha^{-1}\|\varphi\|^2\leq B(\varphi,\varphi),\quad
B(\varphi,\psi)\leq \alpha\|\varphi\|\|\psi\|
\end{equation}
for some positive constant $\alpha>0$. Let us denote
$(\cdot,\cdot)_2 $ the canonical inner product on $L^2(\Omega;\P) $.

From now on, we will say that the unbounded operator $L$ on
$L^2(\Omega;\P)$ is constructed from
$(\Omega,\P,\left<\cdot,\cdot\right>,B)$ if it is constructed as
follows. We consider the inner product $\Pi $ on ${\cal D}\times
{\cal D}$ defined by
$$\Pi(\varphi,\psi)=(\varphi,\psi)_2+\left<\varphi,\psi\right> $$
and the closure $\H $ of ${\cal D} $ with respect to the
corresponding norm. For each $\lambda>0$, the bilinear form
$B_\lambda$ is defined on $ {\cal D}\times {\cal D}$ by
$$B_\lambda(\varphi,\psi)=\lambda(\varphi,\psi)_2+B(\varphi,\psi). $$
From \eqref{equivB}, $B_\lambda$ obviously extends to $\H\times\H$
(this extension is still denoted by $B_\lambda$). Furthermore, it is
continuous and coercive on $\H\times\H$. Thus it defines a resolvent
operator $G_\lambda:L^2(\Omega,\P)\rightarrow\H$, which is
one-to-one. We can then define $L$ as $\lambda-G_\lambda^{-1} $ with
domain ${\rm Dom}(L)=G_\lambda(L^2(\Omega,\P))$. This definition
does not depend on $\lambda>0$. It is readily seen that a function
$\varphi\in \H$ belongs to ${\rm Dom}(L)$ if and only if the map
$\psi\in \H\mapsto B_\lambda(\varphi,\psi)$ is $L^2(\Omega,\P)$
continuous. In this case, we can find $f\in L^2(\Omega,\P)$ such
that $ B_\lambda(\varphi,\cdot)=(f,\cdot)_2$. Then $L\varphi$
exactly matches $f-\lambda\varphi $. Note that
$B(\varphi,\psi)=-(L\varphi,\psi)_2 $ for any $\varphi\in {\rm
Dom}(L)$ and $\psi\in \H$. We point out that the unbounded operator
$L$ is closed and densely defined. Moreover, its adjoint operator
$L^*$ in $L^2(\Omega;\P) $ coincides with the operator constructed
from $(\Omega,\P,\left<\cdot,\cdot\right>,\check{B})$, where the
bilinear form $\check{B} $ is defined on  $ {\cal D}\times {\cal D}$
by $\check{B}(\varphi,\psi)=B(\psi,\varphi) $. As a consequence
$(L^*)^*=L$.

\vspace{2mm}\textit{\textbf{Notations.} In what follows, the
notation $(\H,L,{\rm
Dom}(L),(G_\lambda)_{\lambda>0})=\Xi((\Omega,\P,\left<\cdot,\cdot\right>,B))
$ means that $\H$, $L$, ${\rm Dom}(L)$, $(G_\lambda)_{\lambda>0}$
are constructed from $(\Omega,\P,\left<\cdot,\cdot\right>,B)$ as
explained above.}

%%%%%%%%%%%%%%%%%%%%%%%%%%%%%%%%%%%%%%%%%%%%%%%%%%%%%%%%%%%%%%%%%%%%%%%%%%%%%%%%%%%%%%%%%%%%%%%%%%%%%%%%

\section{Auxiliary Problems}\label{sec:auxiliary}

%%%%%%%%%%%%%%%%%%%%%%%%%%%%%%%%%%%%%%%%%%%%%%%%%%%%%%%%%%%%%%%%%%%%%%%%%%%%%%%%%%%%%%%%%%%%%%%%%%%%%%%%%%%
%\subsection*{Setup}
%%%%%%%%%%%%%%%%%%%%%%%%%%%%%%%%%

\textit{\textbf{Setup and notations.}} Let us now focus on the
different operators induced
 on the random medium $\Omega$ by the matrices ${\bfm a}(\cdot,y)$ and ${\bfm H}(\cdot,y)$, for each $y \in
 \reel^d$. We aim at extending the following operators defined on ${\cal
 C}$ by
\begin{equation}\label{def:oper}
{\bfm S}^y \equiv \frac{1}{2}\sum_{i,j=1}^d D_i\big({\bfm
a}_{ij}(\cdot,y)D_j\,\big), \quad {\bfm L}^y \equiv
\frac{1}{2}\sum_{i,j=1}^d D_i\big(({\bfm a}+{\bfm
H})_{ij}(.,y)D_j\,\big),
\end{equation}
according to the method detailed in Section \ref{sec:construction}.

The positive symmetric bilinear form $(\cdot,\cdot)_1$ is defined on
${\cal C}\times {\cal C}$ by
\begin{equation}
\begin{split}
%({\bfm \varphi} ,{\bfm \psi})_{1,y} & \equiv -({\bfm \varphi} ,{\bfm
%S}(\cdot,y){\bfm \psi})_{2} =(1/2) \bigl(  {\bfm a}(\cdot,y) D{\bfm
%\varphi} , D{\bfm \psi} \bigr)_{2},
%\\
({\bfm \varphi}, {\bfm \psi})_1 & \equiv - ({\bfm \varphi},
\tilde{\bfm S} {\bfm \psi})_2 = (1/2) \bigl(\tilde{\bfm a} D{\bfm
\varphi}, D{\bfm \psi}  \bigr)_{2},
\end{split}
\end{equation}
and the associated seminorm $\|\cdot\|_1 $ by
%$\|{\bfm\varphi}\|_{1}^2
%\equiv ({\bfm \varphi},{\bfm \varphi})_{1,y}$ and
$\|{\bfm \varphi}\|_{1}^2 \equiv ({\bfm \varphi},{\bfm
\varphi})_{1}$.

For any ${\bfm \varphi} ,{\bfm \psi}\in  {\cal C}$, we define the
bilinear forms ($y$ is fixed)
\begin{align*}
{\mathcal B}^{{\small S}} ({\bfm \varphi} ,{\bfm \psi}) & \equiv
-({\bfm S}^y{\bfm \varphi} ,{\bfm \psi})_2=(1/2)({\bfm a}
(\cdot,y)D{\bfm \varphi},D{\bfm \psi})_2,\\
 {\mathcal B}^{{\small L}} ({\bfm \varphi} ,{\bfm \psi}) & \equiv
-({\bfm L}^y{\bfm \varphi} ,{\bfm \psi})_2=(1/2)(({\bfm a}+{\bfm H})
(\cdot,y)D{\bfm \varphi},D{\bfm \psi})_2.
\end{align*}
From Assumption \ref{hyp:cont} and the antisymmetry of ${\bfm H}$,
it is readily seen that $M^{-1} \|{\bfm \varphi}\|_{1}^2\leq
{\mathcal B}^{{\small S}} ({\bfm \varphi} ,{\bfm \varphi}) $ (resp.
$M^{-1} \|{\bfm \varphi}\|_{1}^2\leq {\mathcal B}^{{\small L}}
({\bfm \varphi} ,{\bfm \varphi}) $) and ${\mathcal B}^{{\small S}}
({\bfm \varphi} ,{\bfm \psi})\leq M\|{\bfm \varphi}\|_1\|{\bfm
\psi}\|_1$ (resp. ${\mathcal B}^{{\small L}} ({\bfm \varphi} ,{\bfm
\psi})\leq 2M\|{\bfm \varphi}\|_1\|{\bfm \psi}\|_1$). We can then
define
\begin{align*}
(\H_1,{\bfm S}^y,{\rm Dom}({\bfm
S}^y),(G^{S^y}_\lambda)_{\lambda>0})=
& \Xi(\Omega,\mu,(\cdot,\cdot)_1,{\mathcal B}^{{\small S}}),\\
(\H_1,{\bfm L}^y,{\rm Dom}({\bfm
L}^y),(G^{L^y}_\lambda)_{\lambda>0})= &
\Xi(\Omega,\mu,(\cdot,\cdot)_1,{\mathcal B}^{{\small L}}).
\end{align*}
 Let us additionally denote by $({\bfm
L}^y)^*$ the adjoint operator of ${\bfm L}^y $ in $L^2(\Omega)$.
Note that ${\bfm S}^y$ is self-adjoint.

We define the space $\D$ as the closure in $(L^2(\Omega))^d$ of the
set $\{\tilde{\bfm \sigma }^*D{\bfm \varphi };{\bfm \varphi }\in
{\cal C}\}$. We point out that, whenever $ {\bfm \varphi },{\bfm
\psi }$ belong to $ {\cal C}$,  $2({\bfm \varphi },{\bfm \psi
})_1=(\tilde{\bfm \sigma }^*D{\bfm \varphi },\tilde{\bfm \sigma
}^*D{\bfm \psi  })_2,$ so that the application $\Theta : {\cal C}
\rightarrow  \D, \ {\bfm \varphi } \mapsto  \tilde {\bfm \sigma
}^*D{\bfm \varphi }$ can be extended to the whole space $\H_1$. For
each function ${\bfm f}\in \H_1$, we will note $\nabla^{\tilde
{\sigma}} {\bfm f}$ for $\Theta ({\bfm f})$ and this represents in a
way the gradient of the function ${\bfm f}$ along the direction
$\tilde{\bfm \sigma }$. Similarly, for each fixed $y\in\reel^d$, we
define for any ${\bfm \varphi}\in\H_1$ the gradient along the
direction ${\bfm \sigma}(\cdot,y)$. It will be denoted by
$\nabla^{\sigma(.,y)}{\bfm \varphi}$ and is equal to ${\bfm
\sigma}(\cdot,y)^*D{\bfm \varphi}$ for any ${\bfm \varphi}\in {\cal
C}$. From Assumption \ref{hyp:cont}, for each ${\bfm \varphi}\in
\H_1 $, the mapping $y\in\R^d \mapsto \nabla^{\sigma(.,y)}{\bfm
\varphi}\in \D $ is continuous:
\begin{equation}
\forall(y,h)\in(\reel^d)^2,\quad |\nabla^{\sigma(.,y+h)}{\bfm
\varphi}-\nabla^{\sigma(.,y)}{\bfm \varphi}|_2^2\leq M|h|^2\|{\bfm
\varphi}\|_1^2.
\end{equation}

For $y\in\reel^d$ and ${\bfm \varphi },{\bfm \psi }\in {\cal C}$, we
derive from Assumption \ref{hyp:cont}
\begin{equation}
\label{ty} ( {\bfm L}^y{\bfm \varphi } , {\bfm \psi } )_{2} =
-\frac{1}{2} \bigl( D{\bfm \varphi},({\bfm a} + {\bfm H})(\cdot ,y)
D{\bfm \psi } \bigr) _2\leq C | \nabla^{\tilde{\sigma}} {\bfm
\varphi}|_2 |\nabla^{\tilde{\sigma}}{\bfm \psi }|_2,
\end{equation}
so that we can define a bilinear form ${\bfm T}^y$ on the whole
space $\D \times \D$ such that $\forall {\bfm \varphi },{\bfm \psi
}\in {\cal C}$
\begin{equation}
\label{ty2} - ( {\bfm L}^y{\bfm \varphi } , {\bfm \psi } )_{2} =
{\bfm T}^y(\nabla^{\tilde \sigma} {\bfm \varphi },\nabla^{\tilde
\sigma} {\bfm \psi }).
\end{equation}
Thanks to Assumption \ref{hyp:reg}, we can consider the differential
$\partial {\bfm T}^y$ of ${\bfm T}^y$ defined, for ${\bfm \varphi}, {\bfm
\psi} \in {\mathcal C}$, by $\partial {\bfm T}^y({\bfm \varphi},{\bfm
\psi}) =  \partial_y ({\bfm T}^y({\bfm \varphi},{\bfm \psi}))$. From
Assumption \ref{hyp:cont} and similarly to \eqref{ty}, $\partial {\bfm T}^y$ extends to
${\mathbb D} \times {\mathbb D}$. From Assumption \ref{hyp:cont}, it is then plain to see that the relation $\partial_y\big({\bfm T}^y({\bfm \xi},{\bfm \zeta})\big)= \partial{\bfm T}^y({\bfm \xi},{\bfm \zeta})$ still holds for
${\bfm \xi}, {\bfm \zeta} \in \D$.

Whenever a function ${\bfm b}$ satisfies the property:
\begin{equation}
\exists C>0,\forall{\bfm \varphi}\in {\cal C}, \quad ( {\bfm
b},{\bfm \varphi})_2\leq C\| {\bfm \varphi}\|_1,
\end{equation}
we will say that ${\bfm b}\in \H_{-1}$ and we will define $\|{\bfm
b}\|_{-1}$ as the smallest constant $C$ satisfying this property.

\vspace{2mm} \textbf{\textit{Solvability and regularity of the
resolvent equation.}}  For ${\bfm h}\in L^2(\Omega)$, ${\bfm
u}_{\lambda }(.,y) \equiv G^{L^y}_{\lambda} {\bfm h}$ belongs to
$\H_1\cap{\rm Dom}({\bfm L}^y)$ and satisfies $ \lambda {\bfm
u}_{\lambda }(\cdot,y)-{\bfm L}^y{\bfm u}_{\lambda }(\cdot,y)={\bfm
h}$. Suppose that the right-hand side ${\bfm h}={\bfm h}(\cdot,y)$
depends on the parameter $y\in\R^d$. We now investigate the
$y$-regularity of ${\bfm u}^{\lambda}(\cdot,y)$ from the regularity
of $y \mapsto {\bfm h}(\cdot,y)$ with respect to the norms
$|\cdot|_2$ and $\|\cdot\|_{-1}$. We claim

\begin{proposition}\label{propregloc}
Let us consider ${\bfm h}:y\in \reel^d\mapsto {\bfm h}(.,y)\in
L^2(\Omega)$ and ${\bfm f}:y\in \reel^d\mapsto {\bfm f}(.,y)\in
 L^2(\Omega)\cap\H_{-1}$. Suppose that there exist $C_2,C_{-1}$ such that:

1) the application $y\mapsto{\bfm h}(.,y)\in L^2(\Omega)$ is two
times continuously differentiable in $L^2(\Omega)$. The derivatives
up to order 2 are bounded by $C_2$ in $L^2(\Omega)$ and are
$C_2$-Lipschitz in $L^2(\Omega)$.

2) the application $y\mapsto {\bfm f}(.,y)\in
L^2(\Omega)\cap\H_{-1}$ is two times continuously differentiable in
$\H_{-1}$. The derivatives up to order 2 are bounded by $C_{-1}$ in
$\H_{-1}$ and are $C_{-1}$-Lipschitz in $\H_{-1}$.

Then, for any $\lambda>0 $, the solution ${\bfm u}_{\lambda
}(.,y)\in\H_1\cap {\rm Dom }({\bfm L}^y)$ of the equation
\begin{equation}\label{eqy}
\lambda {\bfm u}_{\lambda }(.,y)-{\bfm L}^y{\bfm u}_{\lambda
}(.,y)={\bfm h}(.,y)+{\bfm f}(.,y)
\end{equation}
is two times continuously differentiable in $\H_1$ with respect to
the parameter $y\in\reel^d$. Furthermore there exists a constant
$D_{\ref{propregloc}}>0$, which only depends on $M,C_{-1}$, such
that the functions ${\bfm g}_\lambda(.,y)={\bfm u}_{\lambda }(.,y)
$, $\partial _y{\bfm u}_{\lambda }(.,y) $, $\partial^2 _{yy} {\bfm
u}_{\lambda }(.,y) $ satisfy the property: $\forall (y,h)\in
\reel^2,$
\begin{subequations}
\begin{align}
\lambda|{\bfm g}_{\lambda }(.,y)|^2_{2}+\|{\bfm g}_{\lambda
}(.,y)\|^2_{1} & \leq
D_{\ref{propregloc}}(1+C^2_2/\lambda),\label{kolmol1}\\
\lambda|{\bfm g}_{\lambda }(.,y+h)-{\bfm g}_{\lambda
}(.,y)|_2^2+\|{\bfm g}_{\lambda }(.,y+h)-{\bfm g}_{\lambda
}(.,y)\|_{1}^2 & \leq
D_{\ref{propregloc}}(1+C^2_2/\lambda)|h|^2\label{kolmol2}.
\end{align}
\end{subequations}
\end{proposition}

\vspace{2mm} \noindent {\bf Proof:} The proof is readily adapted
from \cite[Prop. 4.1]{rhodes:06}. The method consists in
differentiating the resolvent equation (\ref{eqy}) with respect to
the parameter $y\in\R^d$. In the uniformly elliptic setup
\cite[Prop. 4.1]{rhodes:06}, this can be carried out thanks to the
differentiability and the boundedness of ${\bfm a},{\bfm H}$ and
their derivatives up to order 2. In the degenerate setup, we need to
control the matrices ${\bfm a}$ and ${\bfm H}$, as well as their
derivatives up to order $2$ with respect to the parameter
$y\in\R^d$, by the matrix $\widetilde{{\bfm a}} $ (see Assumption
\ref{hyp:cont}) in order to differentiate the function $y\mapsto
{\bfm u}_\lambda(\cdot,y)$ in $\H_1$. \qed

\vspace{2mm} \textbf{\textit{Auxiliary problems: construction of the
correctors.}} The end of this section is now devoted to the study of
the solutions of the so-called auxiliary problems, that means the
solutions ${\bfm u}^i_\lambda (.,y)$ $(i=1,\dots,d)$ of the
resolvent equations
\begin{equation}\label{def:corr}
\lambda {\bfm u}^i_\lambda (.,y)-{\bfm L}^y{\bfm u}^i_\lambda
(.,y)={\bfm b}_i(.,y),
\end{equation}
 where ${\bfm b}_i(.,y)=(1/2)\sum_{j=1}^dD_j\big[({\bfm a}+{\bfm
 H})_{ji}(.,y)\big]$.
The weak form of the resolvent equation then reads for ${\bfm
\varphi}\in {\cal C}$
\begin{align}\label{weakresolvent}
\lambda({\bfm u}^i_\lambda(.,y),{\bfm \varphi})_2 & +{\bfm
T}^y(\nabla^{\widetilde{\sigma}}{\bfm
u}^i_\lambda(.,y),\nabla^{\widetilde{\sigma}}{\bfm \varphi}) =
-(1/2)\big(({\bfm a}+{\bfm H})(.,y)e_i,D{\bfm \varphi}\big)_2.
\end{align}
Having in mind to apply Proposition \ref{propregloc}, we first prove
\begin{lemma}
The mapping $y\mapsto {\bfm b}_i(.,y)\in L^2(\Omega)\cap \H_{-1}  $
is two times continuously differentiable in $\H_{-1}$, and the
derivatives are bounded and Lipschitzian in $\H_1$.
\end{lemma}

\vspace{2mm} \noindent {\bf Proof:} First note that for each ${\bfm
\varphi}\in{\cal C}$,
$$({\bfm b}_i(.,y),{\bfm \varphi})_2= -(1/2)\big(({\bfm a}+{\bfm H})(.,y)e_i,D{\bfm \varphi}\big)_2.$$
From Assumption \ref{hyp:cont}, we easily deduce that ${\bfm
b}_i(.,y)\in\H_{-1}$ and that the mapping $y\in\R^d\mapsto {\bfm
b}_i(.,y)\in\H_{-1}$ is bounded and Lipschitzian.

From Assumption \ref{hyp:cont} again, it is readily seen that the
$\H_{-1}$ derivatives of ${\bfm b}_i$ coincide, for $1\leq k \leq
d$, with the classical derivatives $ \partial_{y_k}{\bfm b}_i$ and
$$(\partial_{y_k}{\bfm b}_i(.,y),{\bfm \varphi})_2= -(1/2)\big((\partial_{y_k}{\bfm a}+
\partial_{y_k}{\bfm H})(.,y)e_i,D{\bfm \varphi}\big)_2\leq C \|{\bfm \varphi}\|_1.$$
Since $ \partial_{y_k}{\bfm a}(\omega)$ and $ \partial_{y_k}{\bfm
H}(\omega)$ are $(M,\tilde{{\bfm a}}(\omega)) $-controlled, the
derivatives are bounded and Lipschitzian in $ \H_1$. The same job
can be carried out for the second order derivatives. Details are
left to the reader.\qed

From Proposition \ref{propregloc} (with ${\bfm h}=0$ and ${\bfm
f}={\bfm b}_i$), the mapping $y\mapsto {\bfm u}^i_\lambda(.,y) $ is
two times continuously differentiable in $\H_{1}$. We now
investigate the asymptotic behavior of ${\bfm u}^i_\lambda$ as well
as its derivatives, as $\lambda$ goes to zero.
\begin{proposition}\label{correctorlocal}
For each fixed $y\in \reel^d$ and $1\leq i \leq d $, the family
$(\nabla^{\widetilde{\sigma}} {\bfm u}_\lambda^i(.,y))_\lambda$
converges to a limit $\widetilde{{\bfm \xi}}_i(.,y)\in
L^2(\Omega)^d$ as $\lambda $ goes to $0$. The same property holds
for the derivatives, namely that the families
$(\nabla^{\widetilde{\sigma}}\partial_{y_j}{\bfm
u}^i_\lambda)_\lambda$, $
(\nabla^{\widetilde{\sigma}}\partial^2_{y_jy_k}{\bfm u}^i_\lambda)_\lambda$
($1\leq i,j,k \leq d$) respectively converge to
$\partial_{y_j}\widetilde{{\bfm \xi}}_i(.,y)$,
$\partial^2_{y_jy_jk}\widetilde{{\bfm \xi}}_i(.,y)$ in
$L^2(\Omega)^d$. Furthermore, we have
\begin{equation*}\label{eq:convcorr}
  \lambda|{\bfm u}_\lambda^i(.,y)|_2^2+\lambda|\partial_{y_j}{\bfm u}_\lambda^i(.,y)|_2^2+\lambda|\partial^2_{y_jy_k}{\bfm u}_\lambda^i(.,y)|_2^2
  \rightarrow 0, \quad \text{ as }\lambda\text{ tends to }0,
\end{equation*}
and, each function ${\bfm g}_\lambda(.,y)={\bfm
u}_\lambda^i(.,y),\partial_{y_k}{\bfm
u}_\lambda^i(.,y),\partial_{y_ky_l}{\bfm u}_\lambda^i(.,y)$
satisfies the property:
\begin{align}
\lambda|{\bfm g}_\lambda(.,y)|_2^2+\|{\bfm g}_\lambda(.,y)\|_1^2 &
\leq C_{\ref{correctorlocal}}\label{eq:borncorr1}\\ \lambda|{\bfm
g}_\lambda(.,y+h)-{\bfm g}_\lambda(.,y)|_2^2+\|{\bfm
g}_\lambda(.,y+h)-{\bfm g}_\lambda(.,y)\|_1^2 & \leq
C_{\ref{correctorlocal}}|h|^2\label{eq:borncorr2}
\end{align}
 for every $ y,h\in
\reel^d$, where $C_{\ref{correctorlocal}}$ is a positive constant
independent of $\lambda>0$ and $y\in\reel^d$.
\end{proposition}

\vspace{2mm} \noindent {\bf Proof:}  The proof does not deeply
differ from Proposition 4.3 in \cite{rhodes:06}, but we nevertheless
set it out because of its importance. From \eqref{kolmol1} (note
that $C_2=0$), we get $\lambda|{\bfm u}^i_\lambda
(.,y)|_2^2+|\nabla^{\widetilde{\sigma}}{\bfm u}^i_\lambda
(.,y)|_2^2\leq C$. Denote by $\widetilde{{\bfm \xi}}_i(.,y)\in
L^2(\Omega)^d$ a weak limit of the family
$(\nabla^{\widetilde{\sigma}}{\bfm u}^i_\lambda (.,y))_\lambda$ as
$\lambda$ goes to $0$. Passing to the limit in
(\ref{weakresolvent}), it is plain to see that $\forall {\bfm
\varphi}\in {\cal C}$
\begin{equation}\label{limcorrec}
 {\bfm T}^y(\widetilde{{\bfm \xi}}_i(.,y),\nabla^{\widetilde{\sigma}}{\bfm
 \varphi})=-(1/2)\big(({\bfm a}+{\bfm H})(.,y)e_i,D{\bfm
 \varphi}\big)_2.
\end{equation}
Since ${\bfm T}^y $ is coercive on $\D\times\D$, this proves the
uniqueness of the weak limit in $\D$. Gathering
(\ref{weakresolvent}) and (\ref{limcorrec}), we get
\begin{equation}\label{cvstrong}
\lambda({\bfm u}^i_\lambda(.,y),{\bfm \varphi})_2  +{\bfm
T}^y(\nabla^{\widetilde{\sigma}}{\bfm
u}^i_\lambda(.,y),\nabla^{\widetilde{\sigma}}{\bfm \varphi}) ={\bfm
T}^y(\widetilde{{\bfm \xi}}_i(.,y),\nabla^{\widetilde{\sigma}}{\bfm
 \varphi}).
 \end{equation} 
 Choosing ${\bfm u}^i_\lambda(.,y)= {\bfm
 \varphi}$ yields:
$$\lambda|{\bfm u}^i_\lambda
(.,y)|_2^2+{\bfm T}^y\big(\nabla^{\widetilde{\sigma}}{\bfm
u}^i_\lambda (.,y),\nabla^{\widetilde{\sigma}}{\bfm u}^i_\lambda
(.,y)\big)\leq {\bfm T}^y\big(\widetilde{{\bfm
\xi}}_i(.,y),\widetilde{{\bfm \xi}}_i(.,y)\big)+\epsilon(\lambda),$$
where the function $\epsilon(\lambda)$ exactly matches ${\bfm
T}^y\big(\widetilde{{\bfm
\xi}}_i(.,y),\nabla^{\widetilde{\sigma}}{\bfm
u}^i_\lambda(.,y)-\widetilde{{\bfm \xi}}_i(.,y)\big)$ and thus
converges to $0$ as $\lambda$ goes to $0$. Hence
$\limsup_{\lambda\rightarrow 0}{\bfm
T}^y\big(\nabla^{\widetilde{\sigma}}{\bfm u}^i_\lambda
(.,y),\nabla^{\widetilde{\sigma}}{\bfm u}^i_\lambda (.,y)\big)\leq
{\bfm T}^y\big(\widetilde{{\bfm \xi}}_i(.,y),\widetilde{{\bfm
\xi}}_i(.,y)\big)$. Denote by ${\bfm T}^S$ the symmetric part of ${\bfm T}^y$
$${\bfm T}^S({\bfm \varphi},{\bfm \psi})=(1/2)\big[{\bfm T}^y({\bfm \varphi},{\bfm \psi})+{\bfm T}^y({\bfm \psi},{\bfm \varphi})\big],\quad {\bfm \varphi},{\bfm \psi}\in\D.$$
From Assumption \ref{hyp:cont} and the antisymmetry of ${\bfm H}$, we have
$$M^{-1}(\widetilde{\bfm \sigma}^*D{\bfm \varphi}, \widetilde{\bfm \sigma}^*D{\bfm \varphi})_2\leq {\bfm T}^S (\nabla^{\widetilde{\bfm \sigma}}{\bfm \varphi},\nabla^{\widetilde{\bfm \sigma}}{\bfm \varphi})\leq M(\widetilde{\bfm \sigma}^*D{\bfm \varphi}, \widetilde{\bfm \sigma}^*D{\bfm \varphi})_2,\quad {\bfm \varphi}\in{\cal C}.$$ 
By density arguments, the quadratic form associated to ${\bfm T}^S $ defines a norm on $\D$ equivalent to the canonical inner product. Moreover, we have just proved that the family $(\nabla^{\widetilde{\sigma}}{\bfm u}^i_\lambda (.,y))_\lambda$ is weakly convergent in $\D$ to $ \widetilde{{\bfm \xi}}_i(.,y)$ and $\limsup_{\lambda\rightarrow 0}{\bfm
T}^S\big(\nabla^{\widetilde{\sigma}}{\bfm u}^i_\lambda
(.,y),\nabla^{\widetilde{\sigma}}{\bfm u}^i_\lambda (.,y)\big)\leq
{\bfm T}^S\big(\widetilde{{\bfm \xi}}_i(.,y),\widetilde{{\bfm
\xi}}_i(.,y)\big)$. Thus the convergence is strong with respect to the norm on $\D$ associated to ${\bfm T}^S$, and consequently  $(\nabla^{\widetilde{\sigma}}{\bfm u}^i_\lambda (.,y))_\lambda$ strongly converges in $(L^2(\Omega))^d$ to $ \widetilde{{\bfm \xi}}_i(.,y)$. From this together with \eqref{cvstrong}, we get
$$ \lambda|{\bfm u}^i_\lambda
(.,y)|_2^2+|\nabla^{\widetilde{\sigma}}{\bfm u}^i_\lambda
(.,y)-\widetilde{{\bfm \xi}}_i(.,y)|_2^2\rightarrow 0\text{ as
}\lambda\rightarrow 0.$$ This proves the first part of the statement
for the function ${\bfm u}^i_\lambda (.,y)$. The second part results
from Proposition \ref{propregloc}, statements \eqref{kolmol1} and
\eqref{kolmol2} (with $C_2=0$). The same job can be carried out for
the successive derivatives of ${\bfm u}^i_\lambda (.,y)$ up to order
$2$. \qed
%%%%%%%%%%%%%%%%%%%%%%%%%%%%%%%%%%%%%%%%%%%%%%%%%%%%%%%%%%%%%%%%%%%%%%%%%%%%%%%%%%%%%%%%%%%%%%%%%%%%%%%%

\section{Dynamics of the process $X^\varepsilon$. Preliminary results}\label{sec:prelim}

%%%%%%%%%%%%%%%%%%%%%%%%%%%%%%%%%%%%%%%%%%%%%%%%%%%%%%%%%%%%%%%%%%%%%%%%%%%%%%%%%%%%%%%%%%%%%%%%%%%%%%%%%%%
\textbf{\textit{Notations.}} \textit{All the results of this section
are valid for any value of the parameter $\varepsilon$. However, to
simplify the notations, we choose $\varepsilon=1$ and thus remove
the parameter $\varepsilon$ from the notations. So the process $X$
stands for the process $X^\varepsilon$ defined by \eqref{eq:diff}.
Finally we denote by $\P_V$ the probability measure $
e^{-2V(y)}\,dy\otimes d\mu$ on $\Omega\times\R^d $ and by $\M_V$ the
coresponding expectation. }

\vspace{2mm} This section is devoted to the study of the
$\Omega\times\R^d $-valued process $(\tau_{X}\omega,X)$, such as its
invariant distribution and the Itô formula. Since these properties are more easily established when the process $X$ possesses regularizing properties, namely that the diffusion coefficient ${\bfm a}$ is uniformly elliptic, most of the following proofs are carried out through
vanishing viscosity methods, that is, in considering a family of
non-degenerate diffusion processes that converges to $X$.

\vspace{2mm} \textit{\textbf{Invariant distribution.}} Let us
introduce a standard d-dimensional Brownian motion $\tilde{B}$
independent of $B$. For each fixed $(\omega,n)\in\Omega\times
\bar{\nat}^*$ and for any $x\in\R^{d}$, we define the Itô process
$X^n$ as the solution of the SDE (with the convention $n^{-1}=0 $ if
$n=\infty$)
\begin{equation*}%\label{}
  X^{n}_t=x+\int_0^t(b+c-n^{-1}\partial_yV)(\omega,X^{n}_r,X^n_r)\,dr+\int_0^t\sigma(\omega,X^{n}_r,X^n_r)\,dB_r+(n/2)^{-1/2}\tilde{B}_t.
\end{equation*}
Note that, for $n=\infty$, $X^{\infty}$ coincides with the process
$X$. For $n\in \bar{\nat}^*$, the process $X^n$ defines a continuous
semigroup $P^n $ on $C_b(\R^d)$ (continuous bounded functions). Its
generator ${\cal L}^n$ coincides on $C^2(\R^d)$ with
\begin{equation}\label{def:ln}{\cal
L}^n=\frac{1}{2}e^{2V(x)}\sum_{i,j}\partial_{
x_i}\big(e^{-2V(x)}(a+H+n^{-1}{\rm
Id})_{ij}(\omega,x,x)\partial_{x_j}\cdot\big).
\end{equation}
 For $  n\in\nat^*$, it is well-known that the distribution of $X^n_t$ ($t>0$) admits a density
 $p^n(\omega,t,x,\cdot)$ with respect to the Lebesgue measure (cf. \cite[Sect. II.2]{stroock}), which is bounded from above by a
  constant $C$ that only depends on $\Lambda,n,t$.
 Thus the semigroup associated to $X^n$ ($ n\in\nat^*$) continuously extends
 to $L^2(\R^d,e^{-2V(x)}\,dx)$. Let us denote by $({\cal
L}^n)^* $ the adjoint of ${\cal L}^n$ in $L^2(\R^d,e^{-2V(x)}\,dx)$,
which coincides on $C^2(\R^d)$ with
\begin{equation}\label{def:lnet}({\cal
L}^n)^*=\frac{1}{2}e^{2V(x)}\sum_{i,j}\partial_{
x_i}\big(e^{-2V(x)}(a-H+n^{-1}{\rm
Id})_{ij}(\omega,x,x)\partial_{x_j}\cdot\big).
\end{equation}

 Now, for $\varphi,\psi\in
 C^\infty_c(\R^d)$, let us compute $\int_{\R^d}{\cal L}^n
 P^n_t\varphi(x)\psi(x)e^{-2V(x)}\,dx$. From \cite{krylov},
 $P^n_t\varphi\in C^2(\R^d)$ so that ${\cal L}^n
 P^n_t\varphi$ can be computed with the help of \eqref{def:ln}. By
 integrating by parts, we obtain
\begin{equation}\label{eq:symm}
\int_{\R^d}{\cal L}^n
 P^n_t\varphi(x)\psi(x)e^{-2V(x)}\,dx=\int_{\R^d} P^n_t\varphi(x)({\cal L}^n)^*\psi(x)e^{-2V(x)}\,dx.
\end{equation}
Moreover, we have ${\cal L}^n
 P^n_t\varphi=P^n_t{\cal L}^n \varphi\in C_b(\R^d)$.  Choose now a
 function $\varrho\in C^\infty_c(\R^d)$ that matches $1$ over the
 ball $B(0;1)$. Define $\psi_m(x)=\varrho(x/m)$. It is readily seen
 that the sequence $({\cal L}^n\psi_m)_m$ is bounded in
 $L^\infty(\R^d)$ and uniformly converges to $0$ on the compact
 subsets of $\R^d$. Thus, choosing $\psi=\psi_m$ in
 \eqref{eq:symm}, and passing to the limit as $m$ goes to $\infty$,
 we get
\begin{equation}\label{eq:invar}
\forall \varphi \in C^\infty_c(\R^d),\quad \int_{\R^d}{\cal L}^n
 P^n_t\varphi(x)e^{-2V(x)}\,dx=0.
\end{equation}
In particular, for any $\varphi\in C^\infty_c(\R^d)$,
$\int_{\R^d}P^n_t\varphi(x) e^{-2V(x)}\,dx=\int_{\R^d}\varphi(x)
e^{-2V(x)}\,dx$, in such a way that, by density arguments, the
probability measure $e^{-2V(x)}\,dx$ is invariant for the process
$X^n$ ($n\geq 1$). Then classical arguments of SDE theory ensure
that the sequence of processes $(X^n)_n$ converges in law in
$C([0,T];\R^d)$ to the process $X$ as $n$ goes to $\infty$. We
deduce that $\int_{\R^d}P_t\varphi(x)
e^{-2V(x)}\,dx=\int_{\R^d}\varphi(x) e^{-2V(x)}\,dx$ holds for
$\varphi\in C_b(\R^d)$. The semigroup associated to $X$ thus extends
to $L^p(\R^d;e^{-2V(x)}\,dx)$ for $p\geq 1$ and the probability
measure $e^{-2V(x)}\,dx$ is also invariant for this semigroup.

Finally, for each ${\bfm \varphi}\in C_b(\Omega\times\R^d) $ (i.e.
for each fixed $\omega\in \Omega$, the function $x\mapsto {\bfm
\varphi}(\tau_x\omega,x)$ is continuous and bounded by a constant
independent of $\omega$) and $n\geq 0$, we deduce from the previous
remarks and the invariance of the measure $\mu$ under space
translations that
\begin{equation}\label{eq:inv}
\bar{\E} [{\bfm \varphi}(\tau_{X^n_t}\omega,X^n_t)]=\M_V[{\bfm
\varphi}(\tau_x\omega,x)]=\M_V[{\bfm \varphi}(\omega,x)],
\end{equation}
so that the mapping ${\bfm \varphi}\in C_b(\Omega\times\R^d)\mapsto
P^n_t({\bfm \varphi})=\E_x[{\bfm \varphi}(\tau_{X^n_t}\omega,X^n_t)
]$ continuously extends to $L^p(\Omega\times\R^d;\P_V)$ for any
$p\geq 1$ and \eqref{eq:inv} holds for ${\bfm \varphi}\in
L^p(\Omega\times\R^d;\P_V)$.

\vspace{2mm} \textit{\textbf{Itô's formula.}} We now aim at
establishing the Itô formula to the process $(\tau_X\omega,X)$ and
to the function $(x,y)\mapsto u_\lambda(\omega,x,y)$, where ${\bfm
u}_\lambda$ is the solution of the resolvent equation \eqref{eqy},
with functions ${\bfm h}(.,y)$ and ${\bfm f}(.,y)$ satisfying the
assumptions of Proposition \ref{propregloc}. This latter proposition
describes the regularity of $ u_\lambda$ with respect to the
variable $y$. Due to the possible degeneracies of ${\bfm \sigma}$,
the difficulty actually lies in the regularity with respect to the
parameter $x\in\R^d$. To apply the Itô formula and get round
technical difficulties, we use viscosity methods again, namely that
we look at the operator $\lambda-{\bfm L^y}-n^{-1}{\bfm \Delta}$ for
$n\in\nat^*$. Obviously, there is no difficulty in solving the
corresponding resolvent equation with the techniques used in Section
\ref{sec:auxiliary} (it suffices to replace ${\bfm a}$ by ${\bfm
a}+n^{-1}{\rm Id}$ and to choose $\tilde{\bfm a}={\rm Id} $)
\begin{equation}\label{eq:nresolvent}
\lambda {\bfm u}^{(n)}_\lambda(\cdot,y)- \bigl( {\bfm L}^y+n^{-1}
{\bfm \Delta} \bigr){\bfm u}^{(n)}_\lambda(\cdot,y)={\bfm
h}(\cdot,y)+{\bfm f}(.,y).
\end{equation}

The strategy then consists in applying the Itô formula in the
non-zero viscosity setting and then in letting $n$ tend to $\infty$.
Thanks to the regularizing parameter $n\in\nat^*$ , the Itô formula
holds in the non-zero viscosity setting (cf \cite[Sect.
5]{rhodes:06}). The following formula thus holds
\begin{align}\label{formiton}
d u_\lambda^{(n)}( X^{n }_t , X^{n}_t) \nonumber = & (\lambda
u_\lambda^{(n)}-h-f)(X^{n }_t ,
X^{n}_t)\,dt+[c-n^{-1}\partial_y\overline{V}]\cdot
Du_\lambda^{(n)}(X^{n }_t , X^{n}_t)\,dt\\
&+(\nabla^{\sigma(.,y)}u_\lambda^{(n)})^*(X^{n }_t , X^{n}_t)\,dB_t
+n^{-1/2}(Du_\lambda^{(n)})(X^{n }_t , X^{n}_t)\,d\widetilde{B}_t\\
 &+b\partial_yu_\lambda^{(n)}(X^{n }_t , X^{n}_t)\,dt+
 [c-n^{-1}\partial_yV]\cdot
\partial_yu_\lambda^{(n)}(X^{n }_t , X^{n}_t)\,dt\nonumber\\
&+(\partial_yu_\lambda^{(n)})^*\sigma(X^{n }_t , X^{n}_t)\,dB_t+
n^{-1/2}(\partial_yu_\lambda^{(n)})(X^{n }_t ,
X^{n}_t)\,d\widetilde{B}_t\nonumber\\ & +(1/2){\rm
trace}([a+n^{-1}Id]\partial^2_{yy}u_\lambda^{(n)})(X^{n }_t ,
X^{n}_t)\,dt\nonumber\\ &+{\rm
trace}([a+n^{-1}Id]D\partial_{y}u_\lambda^{(n)})(X^{n }_t ,
X^{n}_t)\,dt\nonumber.
\end{align}
Having in mind to let $n$ tend to $\infty$ in \eqref{formiton}, let
us now describe the behavior of ${\bfm u}^n_\lambda$ as $n$ tends to
$\infty$. We first claim:
\begin{proposition}
\begin{equation}\label{propregn}
\lim_{n\to \infty}\left[|{\bfm u}^{(n)}_\lambda(.,y)-{\bfm
u}_\lambda(.,y)|_2+\|{\bfm u}^{(n)}_\lambda(.,y)-{\bfm
u}_\lambda(.,y)\|_1+n^{-1}|D{\bfm
u}^{(n)}_\lambda(.,y)|_2^2\right]=0,
\end{equation} and that there exists
a constant $D_{\ref{propregnl}}$ (independent of $n$ and
$y\in\reel^d$) such that
\begin{multline}\label{propregnl}
|{\bfm u}^{(n)}_\lambda(\cdot,y+h)-{\bfm
u}^{(n)}_\lambda(\cdot,y)|_2^2+\|{\bfm u}^{(n)}_\lambda(\cdot,y+h)-
{\bfm u}^{(n)}_\lambda(\cdot,y)\|_1^2\\+ n^{-1}|D{\bfm
u}^{(n)}_\lambda(\cdot,y+h)-D{\bfm u}^{(n)}_\lambda(\cdot,y)|_2^2
\leq D_{\ref{propregnl}}|h|^2.
\end{multline}\label{estn}
Moreover, the same properties hold for the sequences
$(\partial_{y_k}{\bfm u}^{(n)}_\lambda)_n,(\partial^2_{y_ky_l}{\bfm
u}^{(n)}_\lambda)_n$ and their corresponding limits
$(\partial_{y_k}{\bfm u}_\lambda)_n,(\partial^2_{y_ky_l}{\bfm
u}_\lambda)_n$, for $1\leq k,l\leq d $. 
\end{proposition}

\noindent {\bf Proof. } Since the proofs of \eqref{propregn} and \eqref{propregnl}  can be adapted from the proof of Proposition \ref{correctorlocal}, we just set out the guiding line of \eqref{propregn}.

To clarify the notations, we forget for a while the dependence on the parameter $y$. First multiply  \eqref{eq:nresolvent} by ${\bfm u}^{(n)}_\lambda$ and integrate with respect to the measure $\mu$ so as to obtain the estimate:
$$\lambda |{\bfm u}^{(n)}_\lambda|^2_2+ |\nabla^{\widetilde{\sigma}}{\bfm u}^{(n)}_\lambda|^2_2+n^{-1}|D{\bfm u}^{(n)}_\lambda|_2^2\leq C$$ for some constant $C$ only depending on $|{\bfm h}|_2^2/\lambda$ and $\|{\bfm f}\|_{-1}^2$. From this estimate, we deduce that the family $(n^{-1}D{\bfm u}^{(n)}_\lambda)_n$ strongly converges to $0$ in $(L^2(\Omega))^d$ as $n\to \infty$ and that, up to extracting a subsequence, the family $({\bfm u}^{(n)}_\lambda)_n $ weakly converges in $\H_1$ as $n\to \infty$. Multiply once again  \eqref{eq:nresolvent} by a test function ${\bfm \varphi}\in {\cal C}$, integrate with respect to the measure $\mu$ and then pass to the limit as $n\to \infty$ to identity the weak limit in $\H_1 $ as being necessarily equal to ${\bfm u}_\lambda $. So the whole family $({\bfm u}^{(n)}_\lambda)_n $ is weakly convergent in $\H_1$  (not up to a subsequence). It just remains to prove that the convergence actually holds in the strong sense. We can integrate \eqref{eq:nresolvent} and \eqref{eqy} against a test function ${\bfm \varphi}\in {\cal C}$. Since the right-hand sides of \eqref{eq:nresolvent} and \eqref{eqy} coincide, this yields:
$$\lambda({\bfm u}^{(n)}_\lambda,{\bfm \varphi})_2+{\bfm T}^y(\nabla^{\widetilde{\sigma}}{\bfm u}^{(n)}_\lambda,\nabla^{\widetilde{\sigma}}{\bfm \varphi})+
n^{-1}(D{\bfm u}^{(n)}_\lambda,D{\bfm \varphi})_2=\lambda ({\bfm u}_\lambda,{\bfm \varphi})_2+  {\bfm T}^y(\nabla^{\widetilde{\sigma}}{\bfm u}_\lambda,\nabla^{\widetilde{\sigma}}{\bfm \varphi}).$$
Choose ${\bfm \varphi}={\bfm u}^{(n)}_\lambda$ and pass to the limit as $n\to \infty$ and get
$$ \lim_{n\to \infty}\Big(\lambda|{\bfm u}^{(n)}_\lambda|_2^2+{\bfm T}^y(\nabla^{\widetilde{\sigma}}{\bfm u}^{(n)}_\lambda,\nabla^{\widetilde{\sigma}}{\bfm u}^{(n)}_\lambda)+
n^{-1}|D{\bfm u}^{(n)}_\lambda|_2^2\Big)=\lambda|{\bfm u}_\lambda|_2^2+{\bfm T}^y(\nabla^{\widetilde{\sigma}}{\bfm u}_\lambda,\nabla^{\widetilde{\sigma}}{\bfm u}_\lambda).$$
As in Proposition \ref{correctorlocal}, this is sufficient to establish the strong convergence of $({\bfm u}^{(n)}_\lambda)_n $ in $\H_1$ and, consequently, the convergence  $n^{-1}|D{\bfm u}^{(n)}_\lambda|_2^2\to 0$ as $n\to \infty$.\qed

We are now in position to conclude.  Going through formula
(\ref{formiton}), we are faced with functionals of type $\int_t^s
g_n(X^{n}_r, X^{n }_r)\,dr $ (concerning the martingale terms, it
suffices to work on their quadratic variations), where $\M_V[|{\bfm
g}_n-{\bfm g}_0|]\rightarrow 0$ as $n$ tends to $\infty$ and
\begin{equation}\label{eq:regn}
 \forall
(y,h)\in\reel^d\times\reel^d,\quad |{\bfm g}_n(.,y+h)-{\bfm
g}_n(.,y)|_2\leq C|h|
\end{equation} where the constant $C$ depends neither on $n\in\nat$ nor $y,h\in\R^d$. From Lemma
\ref{lemma:conv} below, we prove the convergence of the functional
towards $\int_t^s g_0(X_r, X_r)\,dr$ in $\bar{\P}$-probability and
as a consequence the
\begin{theorem}\label{theoremito}
Let ${\bfm h},{\bfm f}$ be two functions satisfying the assumptions
of Proposition \ref{propregloc}. Let ${\bfm u}_\lambda$ be the
solution of the resolvent equation:
$$\lambda{\bfm u}_\lambda(\cdot,y)-{\bfm L}^y{\bfm u}_\lambda(\cdot,y)={\bfm
h}(\cdot,y)+{\bfm f}(\cdot,y).$$
 Then the following Itô formula holds (we reintroduce the parameter $\varepsilon$):
\begin{align*}
\varepsilon du_\lambda(\overline{X}^{\varepsilon}_t, X^{\varepsilon
}_t) = & \varepsilon^{-1}(\lambda
u_\lambda-h-f)(\overline{X}^{\varepsilon }_t, X^{\varepsilon
}_t)\,dt+c\cdot Du_\lambda(\overline{X}^{\varepsilon }_t,
X^{\varepsilon}_t)\,dt\\
& +(\nabla^{\sigma(.,y)}u_\lambda)^*(\overline{X}^{\varepsilon }_t,
X^{\varepsilon}_t)\,dB_t+b\partial_yu_\lambda(\overline{X}^{\varepsilon
}_t, X^{\varepsilon}_t)\,dt\\
&+\varepsilon(\partial_yu_\lambda)^*\sigma(\overline{X}^{\varepsilon}_t,
X^{\varepsilon }_t)\,dB_t+\varepsilon c\cdot
\partial_yu_\lambda(\overline{X}^{\varepsilon }_t,
X^{\varepsilon }_t)\,dt\\ &+(\varepsilon/2){\rm
trace}(a\partial^2_{yy}u_\lambda)(\overline{X}^{\varepsilon}_t,
X^{\varepsilon}_t)\,dt+{\rm
trace}(aD\partial_{y}u_\lambda)(\overline{X}^{\varepsilon }_t,
X^{\varepsilon}_t)\,dt.
\end{align*}
\end{theorem}
\begin{lemma}\label{lemma:conv}
Consider a sequence of functions ${\bfm g}_n\in
L^1(\Omega\times\R^d;\P_V) $ ($n\geq 0$) such that $\M_V[|{\bfm
g}_n-{\bfm g}_0|]\rightarrow 0$ as $n\to\infty$ and for any
$(y,h)\in\reel^d\times\reel^d$, $|{\bfm g}_n(.,y+h)-{\bfm
g}_n(.,y)|_2\leq C|h|$ for some constant $C$ that depends neither on
$n$ nor $y,h\in\R^d$.

Then $\bar{\E}[|g_n(X^{n}_r, X^{n }_r)-g_0(X_r, X_r)|]\to 0 $ as
$n\to 0$.
\end{lemma}

\vspace{1mm}\noindent {\bf Proof:} First, suppose that ${\bfm g}_0$
is bounded. Let us consider a smooth mollifier $p:\R^d\rightarrow
\R$ and $\varrho\in C^\infty_c(\R^d)$ such that $\varrho=1$ over the
ball $B(0;1)$. We define for $m,q\geq 1 $,
$p_m(\cdot)=m^{d}p(m\,\cdot)$, $\varrho_q(\cdot)=\varrho(\cdot/q) $
and ${\bfm g}^{m,q}_0(\omega,x)= \int_{\R^d}{\bfm
g}_0(\tau_{-x'}\omega,x')\varrho_q(x')p_m(x-x')\,dx'$. Then, from
\eqref{eq:inv},
\begin{align*}
\bar{\E}[|g_n(X^{n}_r, X^{n }_r)-g_0(X_r, X_r)|]\leq & \M_V[|{\bfm
g}_n-{\bfm g}_0|]+2\M_V[|{\bfm g}_0^{m,q}-{\bfm g}_0|]\\
& +\bar{\E}[|g^{m,q}_0(X^{n}_r, X^n_r)-g^{m,q}_0(X_r, X_r)|].
\end{align*}
With classical convolution techniques, we can prove that $m,q$ can
be chosen large enough to make the term $2\M_V[|{\bfm
g}_0^{m,q}-{\bfm g}_0|]$ small. Then, from the Lipschitz regularity
of the coefficients (Assumption \ref{hyp:reg}), the classical theory
of SDEs ensures that $\E_x[\sup_{0\leq t \leq T}|X^n_t-X_t|^2]\leq
n^{-1}D $ for some constant $D$ that only depends on $M$, $\Lambda$
and $T$. For each fixed $m,q \geq 1$ and $\omega\in \Omega $, the
function $x\mapsto g^{m,q}_0(x,x)$ is continuous with compact
support so that $\int_{\R^d}\E_x[|g^{m,q}_0(X^{n}_r,
X^n_r)-g^{m,q}_0(X_r, X_r)|]e^{-2V(x)}\,dx\rightarrow 0 $ as $n\to
\infty$. Then, the Lebesgue theorem ( $g^{m,q}_0 $ is bounded
independently from $\omega$ ) proves that
$\bar{\E}[|g^{m,q}_0(X^{n}_r, X^n_r)-g^{m,q}_0(X_r, X_r)|]$
converges to $0$ as $n$ goes to $\infty$. Therefore, $n$ can be
chosen large enough to make this latter term small. Finally, from
the assumptions of the lemma, even if it means considering larger
$n$, the term $\M_V[|{\bfm g}_n-{\bfm g}_0|]$ is small too. The
proof is then easily completed in the case when ${\bfm g}_0$ is
bounded.

 If ${\bfm g}_0$ is not bounded, it suffices to consider for $n\geq 0$ and
$R>0$, ${\bfm g}_n^R=\max(-R;\min({\bfm g}_n;R))$. It is readily
checked that the sequence $ ({\bfm g}_n^R)_n$ still satisfies all
the assumptions of the lemma in such a way that
$\bar{\E}[|g_n^R(X^{n}_r, X^{n }_r)-g^R_0(X_r, X_r)|]\to 0 $ as
$n\to 0$, for each fixed $R>0$. Then, from \eqref{eq:inv}, $
\bar{\E}[|g_n^R(X^{n}_r, X^{n }_r)-g_n(X^{n}_r, X^{n }_r)|]\leq
\M_V[|{\bfm g}_n^R-{\bfm g}_n|]$ and
$$ \lim_{R\to \infty}\lim_{n\to \infty}\M_V[|{\bfm g}_n^R-{\bfm g}_n|]=\lim_{R\to \infty}\M_V[|{\bfm g}_0^R-{\bfm g}_0|]=0.$$
Since we have
\begin{align*}
\bar{\E}[|g_n(X^{n}_r, X^{n }_r)&-g_0(X^{n}_r, X^{n }_r)|]\leq
\bar{\E}[|g_n^R(X^{n}_r, X^{n }_r)-g_n(X^{n}_r, X^{n
}_r)|]\\&+\bar{\E}[|g_n^R(X^{n}_r, X^{n }_r)-g_0^R(X_r,
X_r)|]+\bar{\E}[|g_0^R(X_r, X_r)-g_0(X^{n}_r, X^{n }_r)|],
\end{align*}
the proof is then easily completed in this case too. \qed

%%%%%%%%%%%%%%%%%%%%%%%%%%%%%%%%%%%%%%%%%%%%%%%%%%%%%%%%%%%%%%%%%%%%%%%%%%%%%%%%%%%%%%%%%%%%%%%%%%%%%%%%

\section{Asymptotic Theorems}

%%%%%%%%%%%%%%%%%%%%%%%%%%%%%%%%%%%%%%%%%%%%%%%%%%%%%%%%%%%%%%%%%%%%%%%%%%%%%%%%%%%%%%%%%%%%%%%%%%%%%%%%%%%

\textbf{\textit{Classical ergodic theorem.}} In this section, we aim
at exploiting the asymptotic properties of the process
$X^\varepsilon$, more precisely Assumption \ref{hyp:erg}, in order
to describe the asymptotic behavior of functionals of type
$\int_0^t\Psi(\overline{X}_r^{\varepsilon },X^{\varepsilon}_r)\,dr $
for a suitable locally stationary random field $ \Psi$. The
classical ergodic theory leads us to guess that the local ergodicity
assumption \ref{hyp:erg} makes this functional average with respect
to its first variable. More precisely,
\begin{theorem}
\label{theoremerg1}{\bf (Ergodic Theorem)} Let us consider ${\bfm
\Psi}:\Omega\times\reel^d  \rightarrow \reel$   such that $\M_V[
|{\bfm \Psi}|]<+\infty$. Denoting $\overline{\Psi}(y)=\M [{\bfm
\Psi}(\cdot,y )]$, the following convergence holds:
\begin{equation}
\em\big[\sup_{0\leq s \leq t}|\int_0^s\Psi
(\overline{X}_r^{\varepsilon
},X^{\varepsilon}_r)\,dr-\int_0^s\overline{\Psi}(X^{\varepsilon}_r)\,dr|
^2\big] \xrightarrow[\varepsilon \rightarrow 0]{}0.
\end{equation}
\end{theorem}
\vspace{1mm} \noindent {\bf Proof:} This result can be proved in the
same way as \cite[Th. 6.1]{rhodes:06}. The only difference consists
in establishing: ${\bfm g}\in {\rm Dom}({\bfm L}^{y})\subset \H_1$
and $ {\bfm L}^{y}{\bfm g}=0$ implies that ${\bfm g}$ is constant
$\mu$ almost surely. In the uniformly elliptic setting, it turns out
that the derivatives $D_i{\bfm g}$ reduce to $0$ and, as a
consequence, ${\bfm g}$ is constant. In the degenerate framework, we
need to use Assumption \ref{hyp:erg} as follows. From Assumption
\ref{hyp:cont}, $\|{\bfm g}\|_1^2\leq M\|{\bfm g}\|_{1,y}^2=-({\bfm
g},{\bfm L}^y{\bfm g})_2=0$. In particular, ${\cal B}^{S^y}({\bfm
g},\cdot)=0$. Hence ${\bfm g}\in {\rm Dom}({\bfm S}^y)$ and $ {\bfm
S}^y{\bfm g}=0$. Thus ${\bfm g}$ is constant (Assumption
\ref{hyp:erg}). \qed

\vspace{2mm} \textbf{\textit{Asymptotic theorem for highly
oscillating functionals.}} Theorem \ref{theoremerg1} describes the
asymptotic behavior of functionals of type $\int_0^t\Psi
(\overline{X}_r^{\varepsilon },X^{\varepsilon}_r)\,dr$ in order to
pass to the limit in \eqref{eq:diff}. However, as explained in
\cite{rhodes:06}, additional difficulties arise in the random
setting in comparison with the periodic one. In particular, we must
describe the asymptotic behavior of the functional
$\int_0^t\Psi_\varepsilon(\overline{X}^\varepsilon_r,X^\varepsilon_r)\,dr
$ for a family $({\bfm \Psi}_\varepsilon)_\varepsilon$ that need not
be convergent in $L^1(\Omega\times\R^d;\P_V)$ but satisfies a sort
of uniform Poincaré inequality. Unlike \cite[Theorem
6.3]{rhodes:06}, technical difficulties due to the degeneracy of the
diffusion coefficient ${\bfm a}$ occur. In particular, because of
the lack of Aronson type estimates, the tightness of the process
$X^\varepsilon$ is not obvious. To prove this tightness, all
asymptotic convergences need be established in $C([0,T];\R^d)$ (note
the $\sup $ in \eqref{equ_cv2}). This is one of the main difficulty
of Theorem \ref{theoremerg2} below in comparison with the uniformly
elliptic setting (see \cite[Theorem 6.3]{rhodes:06}). The strategy
consists in expressing
$\int_0^t\Psi_\varepsilon(\overline{X}^\varepsilon_r,X^\varepsilon_r)\,dr
$ as the sum of two martingales thanks to time reversal arguments,
and then in using the Doob inequality. The Poincaré inequality
\eqref{ass:poincare} ensures that the martingales possess suitable
asymptotic properties.

\begin{theorem}{{\bf (Ergodic theorem II)}}
\label{theoremerg2} Let us consider, for each $\varepsilon>0 $, a
function ${\bfm \Psi}_\varepsilon\in L^2(\Omega\times\reel^d;\P_V)$
satisfying the following Poincaré inequality: for any ${\bfm
\varphi}(\omega,x)={\bfm \chi}(\omega)\varrho(x)$, $({\bfm
\chi},\varrho)\in {\cal C}\times
  C^\infty_c(\R^d) $,
\begin{equation}\label{ass:poincare}
 \M_V[{\bfm \Psi}_\varepsilon{\bfm
  \varphi}]\leq C_\varepsilon \big(\M_V[|{\bfm \sigma}^*(D+\varepsilon\partial_y){\bfm
  \varphi}|^2]\big)^{1/2},
\end{equation}
for some family $(C_\varepsilon)_{\varepsilon>0}$ satisfying
$\varepsilon C_\varepsilon\to 0 $ as $\varepsilon\to 0$. Then
\begin{equation}\label{equ_cv2}
\em\big[\sup_{0\leq s \leq t}|\int_0^s\Psi_\varepsilon
(\overline{X}_r^{\varepsilon },X^{\varepsilon}_r)\,dr| ^2\big]
\xrightarrow[\varepsilon \rightarrow 0]{}0.
\end{equation}
\end{theorem}
\vspace{1mm} \noindent {\bf Proof:}  In what follows, we say that
${\bfm \varphi}\in {\cal C}_\Pi $ if ${\bfm \varphi}(\omega,y)={\bfm
\chi}(\omega)\varrho(y)$, where $({\bfm \chi},\varrho)\in {\cal
C}\times  C^\infty_c(\R^d) $. We aim at constructing, as prescribed
in Section \ref{sec:construction}, the unbounded operators on
$L^2(\Omega\times\reel^d;\P_V)$ that coincide on ${\cal C}_\Pi $ for
$n\in\bar{\nat}^* $ with (here we use the convention $n^{-1}=0$ if
$n=\infty$)
\begin{align}
S^{n,\varepsilon}{\bfm \varphi}=& (1/2)e^{2V}\sum_{i,j=1,\dots
,d}(D_i+\varepsilon\partial_{y_i})\big[e^{-2V}({\bfm a}+n^{-1}{\rm
Id})_{ij}(D_j+\varepsilon\partial_{y_j}){\bfm
\varphi}\big],\\
L^{n,\varepsilon}{\bfm \varphi}=&(1/2)e^{2V}\sum_{i,j=1,\dots
,d}(D_i+\varepsilon\partial_{y_i})\big[e^{-2V}({\bfm a}+{\bfm
H}+n^{-1}{\rm Id})_{ij}(D_j+\varepsilon\partial_{y_j}){\bfm
\varphi}\big].
\end{align}
For $\varepsilon>0$, $n\in\bar{\nat}^*$ and $ {\bfm \varphi},{\bfm
\psi}\in {\cal C}_\Pi$, we define the corresponding bilinear forms
\begin{align}\label{scalarh1}
  \left<{\bfm \varphi},{\bfm \psi}\right>_{n,\varepsilon}
& =(1/2)\M_V\big[(D{\bfm \varphi}+\varepsilon\partial_y{\bfm
\varphi})^*({\bfm a}+n^{-1}{\rm Id})(D{\bfm
\psi}+\varepsilon\partial_y{\bfm \psi})\big],\\
\label{scalarh2}
 B_{n,\varepsilon}({\bfm \varphi},{\bfm \psi})
& =(1/2)\M_V\big[(D{\bfm \varphi}+\varepsilon\partial_y{\bfm
\varphi})^*({\bfm a}+{\bfm H}+n^{-1}{\rm Id})(D{\bfm
\psi}+\varepsilon\partial_y{\bfm \psi})\big].
\end{align}
%\begin{equation}\label{scalarh3}
%\Pi^\alpha_{n,\varepsilon}({\bfm \varphi},{\bfm \psi})
%=\alpha\M_V[{\bfm \varphi}{\bfm
% \psi}]+\left<{\bfm \varphi},{\bfm \psi}\right>_{n,\varepsilon}.
% \end{equation}
Clearly, $\left<\cdot,\cdot\right>_{n,\varepsilon}$ is positive
symmetric (denote by $\|\cdot\|_{n,\varepsilon}$ the corresponding
seminorm). Note that, for each fixed $\varepsilon>0$, the seminorms
$(\|\cdot\|_{n,\varepsilon})_{n\in\nat^*}$ are all equivalent.
Moreover, for $n\in\bar{\nat}^* $, $\|{\bfm
\varphi}\|^2_{n,\varepsilon}\leq B_{n,\varepsilon}({\bfm
\varphi},{\bfm \varphi})$ and $B_{n,\varepsilon}({\bfm
\varphi},{\bfm \psi})\leq 2M^2\|{\bfm
\varphi}\|_{n,\varepsilon}\|{\bfm \psi}\|_{n,\varepsilon} $ for any
${\bfm \varphi},{\bfm \psi}\in {\cal C}_\Pi $ (see Assumption
\ref{hyp:cont}). From Section \ref{sec:construction}, we can define
\begin{align*}
(\H_{n,\varepsilon},S^{n,\varepsilon},{\rm
Dom}(S^{n,\varepsilon}),(G_\lambda^{{\small
S},{n,\varepsilon}})_{\lambda>0})&=\Xi(\Omega\times\R^d,\P_V,
\left<\cdot,\cdot\right>_{n,\varepsilon},\left<\cdot,\cdot\right>_{n,\varepsilon}),\\
(\H_{n,\varepsilon},L^{n,\varepsilon},{\rm
Dom}(L^{n,\varepsilon}),(G_\lambda^{{\small
L},{n,\varepsilon}})_{\lambda>0})&=\Xi(\Omega\times\R^d,\P_V,
\left<\cdot,\cdot\right>_{n,\varepsilon},B_{n,\varepsilon}) .
\end{align*}
and we denote by $(L^{n,\varepsilon})^*$ the adjoint operator of
$L^{n,\varepsilon}$ in $L^2(\Omega\times\R^d;\P_V)$.

Let us now consider a family $({\bfm \Psi}_\varepsilon)_\varepsilon$
of functions in $ L^2(\Omega\times\reel^d;\P_V)$ satisfying
\eqref{ass:poincare} for some family
$(C_\varepsilon)_{\varepsilon>0}$ such that $\varepsilon
C_\varepsilon\to 0 $ as $\varepsilon\to 0$. Fix $n\in\bar{\nat}^*$.
Define $ {\bfm \varphi}_{n,\varepsilon}\equiv
G_{\varepsilon^2}^{S,n,\varepsilon}({\bfm
\Psi}_{\varepsilon})$, which satisfies $\varepsilon^2\M_V[{\bfm
\varphi}_{n,\varepsilon}{\bfm \psi}]+ \left<{\bfm
\varphi}_{n,\varepsilon},{\bfm
\psi}\right>_{n,\varepsilon}=\M_V[{\bfm \Psi}_\varepsilon{\bfm
\psi}]$ for any ${\bfm \psi}\in\H_{n,\varepsilon}$. Choosing ${\bfm
\psi}={\bfm \varphi}_{n,\varepsilon}$, using \eqref{ass:poincare} and the standard estimate $ab\leq a^2/2+b^2/2$
leads to 
\begin{align*}%\label{}
 \varepsilon^2\M_V[|{\bfm \varphi}_{n,\varepsilon}|^2]+ \|{\bfm \varphi}_{n,\varepsilon}\|^2_{n,\varepsilon}&=\M_V[{\bfm \Psi}_\varepsilon{\bfm \varphi}_{n,\varepsilon}]\leq C_\varepsilon\sqrt{2} \|{\bfm \varphi}_{n,\varepsilon}\|_{0,\varepsilon}\leq \sqrt{2}  C_\varepsilon \|{\bfm \varphi}_{n,\varepsilon}\|_{n,\varepsilon}\\&\leq C_\varepsilon^2+ \|{\bfm \varphi}_{n,\varepsilon}\|_{n,\varepsilon}^2/2
\end{align*}
in such a way that
\begin{equation}\label{estpoin}
 \varepsilon^2\M_V[|{\bfm \varphi}_{n,\varepsilon}|^2]+ \|{\bfm \varphi}_{n,\varepsilon}\|^2_{n,\varepsilon}/2\leq C_\varepsilon^2.
\end{equation}

Once again, to apply the Itô formula, we use vanishing viscosity
methods in order to get round the lack of regularity of $ {\bfm
\varphi}_{n,\varepsilon}$ because of the degeneracy of ${\bfm a}$.
In the non-degenerate framework ($n\geq 1$), from \cite[Proof of
Lemma 6.3]{rhodes:06}, standard convolution technics provide us with
a $\H_{n,\varepsilon} $-sequence $({\bfm
\varphi}^m_{n,\varepsilon})_{m\in\nat} $ of smooth functions, namely
that for each fixed $\omega\in\Omega$ the function $x\mapsto {\bfm
\varphi}^m_{n,\varepsilon}(\tau_{x/\varepsilon}\omega,x)$ is a
$C^\infty(\R^d)$-function, such that $\M_V[|{\bfm
\varphi}^m_{n,\varepsilon}-{\bfm
\varphi}_{n,\varepsilon}|^2+|S^{n,\varepsilon}{\bfm
\varphi}^m_{n,\varepsilon}-S^{n,\varepsilon}{\bfm
\varphi}_{n,\varepsilon}|^2]+\|{\bfm
\varphi}^m_{n,\varepsilon}-{\bfm
\varphi}_{n,\varepsilon}\|^2_{n,\varepsilon}\rightarrow 0$ as $m$
goes to $\infty$.

We are now going to use a time reversal argument. Let us consider
the process (introduced in Section \ref{sec:prelim})
\begin{equation*}%\label{}
  X^{n,\varepsilon}_t=x+\int_0^t(\varepsilon^{-1}b+c-n^{-1}\partial_yV)(\omega,\overline{X}^{n,\varepsilon}_r,X^{n,\varepsilon}_r)\,dr
  +\int_0^t\sigma(\omega,
  \overline{X}^{n,\varepsilon}_r,X^{n,\varepsilon}_r)\,dB_r+(n/2)^{-1/2}\tilde{B}_t,
\end{equation*}
where $
\overline{X}^{n,\varepsilon}_r=X^{n,\varepsilon}_r/\varepsilon$. As
explained in Section \ref{sec:prelim}, its generator coincides on
$C^2(\R^d)$ with
\begin{equation*}%\label{def:ln}
{\cal L}^{n,\varepsilon}=\frac{e^{2V(x)}}{2}\sum_{i,j}\partial_{
x_i}\big(e^{-2V(x)}(a+H+n^{-1}{\rm
Id})_{ij}(\omega,x/\varepsilon,x)\partial_{x_j}\cdot\big)
\end{equation*}
and admits $e^{-2V(x)}\,dx$ as invariant measure. Furthermore, for a
fixed $T>0$, the generator of the time reversed process $t\mapsto
X^{n,\varepsilon}_{T-t}$ with initial law $e^{-2V(x)}\,dx $
coincides with the adjoint of ${\cal L}^{n,\varepsilon}$ in
$L^2(\R^d;e^{-2V(x)}\,dx) $. For each $\varphi\in C^2(\R^d) $, it
exactly matches
\begin{equation*}%\label{def:ln}
({\cal
L}^{n,\varepsilon})^*\varphi=\frac{e^{2V(x)}}{2}\sum_{i,j}\partial_{
x_i}\big(e^{-2V(x)}(a-H+n^{-1}{\rm
Id})_{ij}(\omega,x/\varepsilon,x)\partial_{x_j}\varphi\big)
\end{equation*}
As a consequence, observe that, for any $ 0\leq s \leq t\leq T$,
\begin{align*}
\varphi^m_{n,\varepsilon}(\overline{X}^{n,\varepsilon}_t,X^{n,\varepsilon}_t)=
&
\varphi^m_{n,\varepsilon}(\overline{X}^{n,\varepsilon}_s,X^{n,\varepsilon}_s)+
\int_s^t[{\cal
L}^{n,\varepsilon}(\varphi^m_{n,\varepsilon}(\cdot/\varepsilon,\cdot))](\overline{X}^{n,\varepsilon}_r,X^{n,\varepsilon}_r)\,dr\\
&+(\overrightarrow{{\cal
M}}^{m,n,\varepsilon}_t-\overrightarrow{{\cal
M}}^{m,n,\varepsilon}_s),
\end{align*}
where $\overrightarrow{{\cal M}}^{m,n,\varepsilon}   $ is a
martingale with respect to the forward filtration $({\cal
F}^{n,\varepsilon} _t)_{0\leq t\leq T}$ and ${\cal
F}^{n,\varepsilon}_t$ is the $\sigma$-algebra on $\reel^d $
generated by $\left\{X^{n,\varepsilon}_r;0\leq r\leq t\right\}$. In
the same way,
\begin{align*}
\varphi^m_{n,\varepsilon}(\overline{X}^{n,\varepsilon}_s,X^{n,\varepsilon}_s)=
&
\varphi^m_{n,\varepsilon}(\overline{X}^{n,\varepsilon}_t,X^{n,\varepsilon}_t)+
\int_s^t[({\cal
L}^{n,\varepsilon})^*(\varphi^m_{n,\varepsilon}(\cdot/\varepsilon,\cdot))](\overline{X}^{n,\varepsilon}_r,X^{n,\varepsilon}_r)\,dr\\&
+(\overleftarrow{{\cal M}}^{m,n,\varepsilon}_t-\overleftarrow{{\cal
M}}^{m,n,\varepsilon}_s),
\end{align*} where $\overleftarrow{{\cal M}}^{m,n,\varepsilon } $ is a martingale with respect
to the backward filtration $({\cal G}^{n,\varepsilon} _t)_{0\leq
t\leq T}$ and ${\cal G}^\varepsilon _s$ is the $\sigma$-algebra on
$\R^d $ generated by $\left\{ X^{n,\varepsilon}_r;t\leq r\leq
T\right\}$. Add these two expressions:
$$-2\varepsilon^{-2}\int_s^tS^{n,\varepsilon}\varphi^m_{n,\varepsilon}(\overline{X}^{n,\varepsilon}_r
,X^{n,\varepsilon}_r)\,dr=(\overrightarrow{{\cal
M}}^{m,n,\varepsilon}_t-\overrightarrow{{\cal
M}}^{m,n,\varepsilon}_s)+(\overleftarrow{{\cal
M}}^{m,n,\varepsilon}_t-\overleftarrow{{\cal
M}}^{m,n,\varepsilon}_s).$$
 We further mention that the quadratic
variations of both martingales exactly match
$$\varepsilon^{-2}\int_s^t[(D+\varepsilon\partial_y) \varphi^m_{n,\varepsilon}]^*a
[(D+\varepsilon\partial_y
)\varphi^m_{n,\varepsilon}\big]^*(\overline{X}^{n,\varepsilon}_r,X^{n,\varepsilon}_r)\,dr,$$
in such a way that the Doob inequality yields $$\em\Big[\sup_{0\leq
s \leq t}|\int_0^sS^{n,\varepsilon}\varphi^m_{n,\varepsilon}
(\overline{X}_r^{n,\varepsilon },X^{n,\varepsilon}_r)\,dr|
^2\Big]\leq 16 T\varepsilon^2\|{\bfm \varphi}^m_{n,\varepsilon}
\|^2_{n,\varepsilon}.$$ Letting $m$ go to $\infty$, reminding
that $\varepsilon^2{\bfm
\varphi}_{n,\varepsilon}-S^{n,\varepsilon}{\bfm
\varphi}_{n,\varepsilon}={\bfm \Psi}_{\varepsilon} $ and using \eqref{estpoin} leads to
$$\em\big[\sup_{0\leq s \leq t}|\int_0^s\Psi_\varepsilon
(\overline{X}_r^{n,\varepsilon },X^{n,\varepsilon}_r)\,dr|
^2\big]\leq 32T \varepsilon^2\|{\bfm \varphi}_{n,\varepsilon}
\|^2_{n,\varepsilon}+2T\varepsilon^4\M_V[|{\bfm
\varphi}_{n,\varepsilon}|^2]\leq 68T \varepsilon^2C_\varepsilon^2.$$
We then complete the proof in letting $n$ go to $\infty$ and in
using the fact that $X^{n,\varepsilon} $ converges in $C([0,T];\R^d)$
towards $X^\varepsilon $ as $n$ goes to $ \infty$.\qed

%%%%%%%%%%%%%%%%%%%%%%%%%%%%%%%%%%%%%%%%%%%%%%%%%%%%%%%%%%%%%%%%%%%%%%%%%%%%%%%%%%%%%%%%%%%%%%%%%%%%%%%%

\section{Proof of Theorem \ref{maintheorem} and Proposition
\ref{propgeom}}\label{sec:proof}

%%%%%%%%%%%%%%%%%%%%%%%%%%%%%%%%%%%%%%%%%%%%%%%%%%%%%%%%%%%%%%%%%%%%%%%%%%%%%%%%%%%%%%%%%%%%%%%%%%%%%

\noindent {\bf Proof of Theorem \ref{maintheorem}.} Section
\ref{sec_tightness} below is devoted to proving the tightness of the
family of processes $(X^\varepsilon)_\varepsilon$ in
$C([0,T];\reel^d)$. It remains to prove that there is a unique
possible weak limit for all converging subsequences.

From now on, the corrector ${\bfm u}_\lambda^i$ ($\lambda>0 $ and
$1\leq i \leq d $) stands for the solution of \eqref{def:corr}.
Applying the Ito formula (Theorem \ref{theoremito}) to the
correctors leads to
\begin{align*}
dX^{\varepsilon}_t= & -\varepsilon
du_{\varepsilon^2}(\overline{X}^{\varepsilon }_t, X^{\varepsilon
}_t)+\varepsilon(\partial_yu_{\varepsilon^2})^*\sigma(\overline{X}^{\varepsilon
}_t, X^{\varepsilon }_t)\,dB_t \\ & +[\varepsilon
u_{\varepsilon^2}+\varepsilon c\cdot
\partial_yu_{\varepsilon^2}+(\varepsilon/2){\rm
trace}(a\partial^2_{yy}u_{\varepsilon^2})](\overline{X}^{\varepsilon
}_t, X^{\varepsilon }_t)\,dt\\
 &+[b\partial_yu_{\varepsilon^2}+c\cdot (I+ Du_{\varepsilon^2})+{\rm
trace}(aD\partial_{y}u_{\varepsilon^2})](\overline{X}^{\varepsilon}_t,
X^{\varepsilon }_t)\,dt\\
&+[\sigma+Du_{\varepsilon^2}\sigma](\overline{X}^{\varepsilon }_t,
X^{\varepsilon }_t)\,dB_t\\ \equiv &
d\Theta_t^{1,\varepsilon}+d\Theta_t^{2,\varepsilon}+d\Theta_t^{3,\varepsilon}+d\Theta_t^{4,\varepsilon}
\end{align*}
Concerning the first term, we have
$\em\big[|\Theta_t^{1,\varepsilon}|^2\big] \leq
(1+T)\varepsilon^2\M_V\big[|{\bfm
u}_{\varepsilon^2}|^2+M^2|\partial_y{\bfm
u}_{\varepsilon^2}|^2\big]$ for $ 0\leq t \leq T$. This latter
quantity converges to $0$ as $\varepsilon$ goes to $0$ from
Proposition \ref{correctorlocal}. The same job can be carried out
for $ \Theta^{2,\varepsilon}$ and the same conclusion holds.

The main difficulty actually lies in the term
$\Theta^{3,\varepsilon}$, especially in the part corresponding to
${\bfm b}\partial_y{\bfm u}_{\varepsilon^2}$. Concerning the
remaining part ${\bfm c}\cdot (I+ D{\bfm u}_{\varepsilon^2})+{\rm
trace}({\bfm a}D\partial_{y}{\bfm u}_{\varepsilon^2}) $, it is
readily seen (see Proposition \ref{correctorlocal}) that it
converges in $L^2(\Omega\times\R^d;\P_V) $ and thus Theorem
\ref{theoremerg1} can be applied. As a consequence, we have
$$\em\big[\sup_{0 \leq t \leq T}\big|\int_0^t[c\cdot (I+ Du_{\varepsilon^2})+{\rm
trace}(aD\partial_{y}u_{\varepsilon^2})](\overline{X}^{\varepsilon}_r,
X^{\varepsilon }_r)\,dr - \int_0^t\bar{\Phi}(X^{\varepsilon
}_r)\,dr\big|^2\big]\rightarrow 0\text{ as }\varepsilon \rightarrow
0,  $$ where $\bar{\Phi}(y)=\lim_{\lambda\rightarrow 0}\M[{\bfm
c}\cdot (I+ D{\bfm u}_{\lambda})+{\rm trace}({\bfm
a}D\partial_{y}{\bfm u}_{\lambda})(\cdot,y) ] $. It remains to treat
the term $({\bfm b}\partial_y{\bfm u}_{\varepsilon^2})_\varepsilon$.
Note that the $L^2$-norm of $ {\bfm b}\partial_y{\bfm
u}_{\varepsilon^2}$ need not be convergent. That is why we have in
mind to use Theorem \ref{theoremerg2}. Up to introducing new
correctors, we will prove that ${\bfm b}\partial_y{\bfm
u}_{\varepsilon^2}$ can be divided into two parts, satisfying
respectively Theorems \ref{theoremerg1} and \ref{theoremerg2}. To
understand how this decomposition occurs, let us consider a test
function ${\bfm \varphi}\in {\cal C}_\Pi$. Then two successive
integrations by parts yield, for $1\leq i,j\leq d$, (we use the
convention of summation over repeated indices)
\begin{eqnarray*}
  \M_V\big[{\bfm b}_j\partial_{y_j}{\bfm u}^i_{\varepsilon^2}{\bfm \varphi}\big] & = & (1/2)\M_V\big[
  D_p({\bfm a}+{\bfm H})_{pj}\partial_{y_j}{\bfm u}^i_{\varepsilon^2}{\bfm \varphi}\big] \\
  & = & -(1/2)\M_V\big[
  ({\bfm a}+{\bfm H})_{pj}(D_p\partial_{y_j}{\bfm u}^i_{\varepsilon^2}{\bfm \varphi}+
  \partial_{y_j}{\bfm u}^i_{\varepsilon^2}D_p{\bfm \varphi})\big] \\
  & = & -(1/2)\M_V\big[
  ({\bfm a}+{\bfm H})_{pj}\big(D_p\partial_{y_j}{\bfm u}^i_{\varepsilon^2}{\bfm \varphi}+
  \partial_{y_j}{\bfm u}^i_{\varepsilon^2}(D_p+\varepsilon\partial_{y_p}){\bfm \varphi}\big)\big] \\
  & & +(\varepsilon/2)\M_V\big[
  ({\bfm a}+{\bfm H})_{pj}\partial_{y_j}{\bfm u}^i_{\varepsilon^2}\partial_{y_p}{\bfm
  \varphi}\big]\\
    & = & -(1/2)\M_V\big[
  ({\bfm a}+{\bfm H})_{pj}D_p\partial_{y_j}{\bfm u}^i_{\varepsilon^2}{\bfm \varphi}\big]-(1/2)\M_V\big[
  \partial_{y_j}{\bfm u}^i_{\varepsilon^2}(D_p+\varepsilon\partial_{y_p}){\bfm \varphi}\big] \\
  & & -(\varepsilon/2)\M_V\big[
  \partial_{y_p}({\bfm a}+{\bfm H})_{pj}\partial_{y_j}{\bfm u}^i_{\varepsilon^2}{\bfm
  \varphi}+({\bfm a}+{\bfm H})_{pj}\partial_{y_jy_p}^2{\bfm u}^i_{\varepsilon^2}{\bfm
  \varphi}\big]\\
  & &+ \varepsilon\M_V\big[({\bfm a}+{\bfm H})_{pj}\partial_{y_j}{\bfm u}^i_{\varepsilon^2}{\bfm
  \varphi}\partial_{y_p}V\big].
\end{eqnarray*}
So, for $1\leq i\leq d$, define the correcting part
$\textrm{Corr}_\varepsilon^i(\omega,y)=(\varepsilon/2)\partial_{y_p}({\bfm a}+{\bfm H})_{pj}\partial_{y_j}{\bfm u}^i_{\varepsilon^2}+(\varepsilon/2)({\bfm a}+{\bfm H})_{pj}\partial_{y_jy_p}^2{\bfm
  u}^i_{\varepsilon^2}-\varepsilon({\bfm a}+{\bfm H})_{pj}\partial_{y_j}{\bfm u}^i_{\varepsilon^2}\partial_{y_p}V$, the $L^2$-converging part $\textrm{Conv}_\varepsilon^i(\omega,y)=
  -(1/2)   ({\bfm a}+{\bfm H})_{pj}D_p\partial_{y_j}{\bfm
  u}^i_{\varepsilon^2}$ and $L^2$-diverging part $\textrm{Div}_\varepsilon^i(\omega,y)=[{\bfm b}_j\partial_{y_j}{\bfm
u}^i_{\varepsilon^2}+\textrm{Corr}_\varepsilon^i-\textrm{Conv}_\varepsilon^i](\omega,y)$.
From the previous calculation, $\textrm{Div}_\varepsilon^i$
satisfies the "Poincaré inequality" \eqref{ass:poincare}, namely
that for any function ${\bfm \varphi}$ in ${\cal C}_\Pi$,
$\M_V\big[{\bfm \Psi}_\varepsilon{\bfm \varphi}\big]\leq
\big(M_V[|\partial_{y}{\bfm
u}^i_{\varepsilon^2}|^2]\big)^{1/2}\big(M_V[|(D+\varepsilon\partial_{y}){\bfm
\varphi}|^2]\big)^{1/2}$. Moreover, Proposition \ref{correctorlocal}
ensures that $\varepsilon\big(M_V[|\partial_{y}{\bfm
u}^i_{\varepsilon^2}|^2]\big)^{1/2}\to $ as $\varepsilon$ goes to
$0$. Consequently, \eqref{equ_cv2} holds for $
\textrm{Div}_\varepsilon^i$. Thanks to Proposition
\ref{correctorlocal}, the family
$(\textrm{Corr}_\varepsilon^i)_\varepsilon$ converges in
$L^2(\Omega\times\R^d;\P_V)$ towards $0$. As a consequence,
$\em\big[\big(\int_0^t\textrm{Corr}_\varepsilon^i(\overline{X}^{\varepsilon
}_r, X^{\varepsilon }_r)\,dr\big)^2\big] $ tends to $0$ as
$\varepsilon$ goes to $0$. Then, Theorem \ref{theoremerg1} ensures
that $\em \big[\sup_{0\leq t \leq T}|\int_0^t
\textrm{Conv}_\varepsilon^i(\overline{X}^\varepsilon,X^\varepsilon_r)\,dr-
  \int_0^t \overline{\Gamma}(X^\varepsilon_r)\,dr|^2\big]\to 0 $ as
  $\varepsilon\to 0$, where $\overline{{\bfm \Gamma}}(y)\equiv \lim_{\lambda\to 0}-(1/2)\M\big[
  ({\bfm a}+{\bfm H})_{pj}D_p\partial_{y_j}{\bfm u}^i_{\lambda}(\cdot,y)\big]$. To sum up, this proves that
\begin{equation}\label{eq_convb}
  \em\big[\sup_{0\leq t \leq T}|\int_0^t b\cdot\partial_{y}u^i_{\varepsilon^2}(\overline{X}^\varepsilon,X^\varepsilon_r)\,dr-
  \int_0^t \overline{\Gamma}(X^\varepsilon_r)\,dr|^2\big]\to 0
\end{equation}
as $\varepsilon$ tends to $0$.

Concerning the martingale part $\Theta^{4,\varepsilon} $, it
suffices to apply Theorem \ref{theoremerg1} to the quadratic
variations.

Hence each possible limit point $X$ in $C[0,t];\reel^d)$ of the
process $X^\varepsilon$ must solve the martingale problem $X_t=
  x+\int_0^t\overline{B}(X_r)\,dr+\int_0^t\overline{A}^{1/2}(X_r)\,dB_r$,
  where the entries of $\bar{B}$ are given by
\begin{align*}
\overline{B}_i &= \lim_{\lambda\to 0}\M\big[-(1/2)({\bfm a}+{\bfm
H})_{pj}D_p\partial_{y_j}{\bfm u}_\lambda^i+{\bfm
c}_j(\delta_{ij}+D_j{\bfm u}_\lambda^i)+{\bfm
a}_{pj}D_j\partial_{y_p}{\bfm u}_\lambda^i\big]\\
&= \lim_{\lambda\to 0}\M\big[(1/2)({\bfm a}+{\bfm
H})_{pj}D_p\partial_{y_j}{\bfm u}_\lambda^i+{\bfm
c}_j(\delta_{ij}+D_j{\bfm u}_\lambda^i)\big]\\
&= \frac{e^{2V}}{2}\partial_{y_j}\big(e^{-2V}\lim_{\lambda\to
0}\M\big[({\bfm a}+{\bfm H})_{pj}(\delta_{ij}+D_p{\bfm
u}_\lambda^i)\big]\big)
\end{align*}

  Thanks to Proposition \ref{propgeom}, it is readily seen that the
coefficients $\overline{B}$ and $\overline{A}^{1/2}$ are two times
continuously differentiable with bounded derivatives up to order
two. In particular, they are Lipschitzian and there exists a unique
solution to the corresponding martingale problem.\qed

\medskip
\noindent {\bf Proof of Proposition \ref{propgeom}.}  The strategy
consists in introducing the homogenized diffusion coefficient
associated to the operator $\widetilde{{\bfm S}}$ and in comparing
it with $\bar{A}(y)$. So we define the $d\times d$ nonnegative
symmetric matrix $\widetilde{A}$ as the unique symmetric matrix
satisfying (this is the classical variational formula for the
homogenized coefficient associated to $\widetilde{{\bfm S}}$, see
\cite{olla:94} for further details)
\begin{equation}\label{varform}
   \forall x\in\R^d,\quad \langle
x,\widetilde{A}x\rangle=\inf_{{\bfm \varphi}\in {\cal
C}}\M\big[|\widetilde{{\bfm \sigma}}^*(D{\bfm \varphi}+x)|^2\big].
\end{equation}
Due to Assumption \ref{hyp:cont}, we have for each function ${\bfm
\varphi}\in {\cal C}$,
\begin{equation*}%\label{}
M^{-1}\langle x,\widetilde{A}x\rangle\leq
M^{-1}\M\big[|\widetilde{{\bfm \sigma}}^*(D{\bfm
\varphi}+x)|^2\big]\leq \M\big[|{\bfm \sigma}^*(\cdot,y)(D{\bfm
\varphi}+x)|^2\big].
\end{equation*}
Since ${\cal C}$ is dense in $\H_1$, we can choose ${\bfm
\varphi}={\bfm u}_\lambda(\cdot,y)\cdot x$ and then pass to the
limit as $\lambda$ tends to $0$. We obtain $M^{-1}\langle
x,\widetilde{A}x\rangle \leq \langle x,\overline{A}(y)x\rangle$.

Now we turn to the auxiliary problems (subsection
\ref{sec:auxiliary}). Denoting by ${\mathbb L}$ the closure of
$\{\tilde{\bfm \sigma} ^*{\bfm \zeta}, \ {\bfm \zeta} \in
L^2(\Omega;\R^d)\}$, we can extend ${\bfm T}^y$ to the whole
${\mathbb L}$ as follows
\begin{equation}\label{def:tyl}
\forall {\bfm \zeta}, {\bfm \theta} \in L^2(\Omega,\R^d), \ {\bfm
T}^y(\tilde{\bfm \sigma}^* {\bfm \zeta}, \tilde{\bfm \sigma}^* {\bfm
\theta}) = (1/2)\bigl([{\bfm a}+{\bfm H}](\cdot,y) {\bfm
\zeta},{\bfm \theta}\bigr)_2.
\end{equation}
The underlying quadratic form is still denoted by ${\bfm
T}^y(\cdot)$. Furthermore, from Assumption \ref{hyp:cont}, for some
positive constant $C$ only depending on $M$, we have
\begin{equation}\label{prop:tyl}
{\bfm T}^y(\tilde{\bfm \sigma}^* {\bfm \zeta}, \tilde{\bfm \sigma}^*
{\bfm \theta})\leq C {\bfm T}^y( \tilde{\bfm \sigma}^* {\bfm
\zeta})^{1/2}{\bfm T}^y(\tilde{\bfm \sigma}^* {\bfm \theta})^{1/2}.
\end{equation}
 Equation (\ref{limcorrec}) then
reads, for any function ${\bfm \varphi}\in{\cal C}$,
\begin{equation}\label{weakty}
\forall x \in \R^d, \ {\bfm T}^y ( \widetilde{{\bfm \xi}} (\cdot,y)
x, \tilde{\bfm \sigma}^* {\bfm D}{\bfm \varphi} ) = -(1/2) ( [{\bfm
a} + {\bfm H}](\cdot,y)x,  {\bfm D} {\varphi} \bigr)_2 = - {\bfm
T}^y( \tilde{\bfm \sigma}^* x, \tilde{\bfm \sigma}^*  {\bfm D}{\bfm
\varphi}).
\end{equation}
From \eqref{hom:a}, \eqref{def:tyl} and \eqref{weakty}, we have for
any function ${\bfm \varphi}\in{\cal C} $
\begin{align*}
\langle x,\overline{A}(y)x\rangle & = 2\lim_{\lambda\to 0}{\bfm T}^y
\big( \widetilde{{\bfm \sigma}}^*x+\nabla^{\widetilde{\sigma}}{\bfm
u}_\lambda(\cdot,y)  x \big)=2{\bfm T}^y ( \widetilde{{\bfm
\sigma}}^*x+\widetilde{{\bfm \xi}}
(\cdot,y) x )\\
& =2{\bfm T}^y ( \widetilde{{\bfm \sigma}}^*x+\widetilde{{\bfm \xi}}
(\cdot,y) x, \widetilde{{\bfm \sigma}}^*x+\widetilde{\bfm \sigma}^*
{\bfm D}{\bfm \varphi} )\\ &  \leq 2 C {\bfm T}^y ( \widetilde{{\bfm
\sigma}}^*x+\widetilde{{\bfm \xi}} (\cdot,y) x )^{1/2}{\bfm T}^y (
\widetilde{{\bfm \sigma}}^*x+\widetilde{\bfm \sigma}^* {\bfm D}{\bfm
\varphi})^{1/2}.
\end{align*}
Gathering this with the inequality ${\bfm T}^y ( \widetilde{{\bfm
\sigma}}^*x+\widetilde{\bfm \sigma}^* {\bfm D}{\bfm \varphi} )\leq M
\M\big[|\widetilde{{\bfm \sigma}}^*x+\widetilde{\bfm \sigma}^* {\bfm
D}{\bfm \varphi}|^2\big]$ and \eqref{varform}, we deduce $\langle
x,\overline{A}(y)x\rangle\leq  2C^2M\langle x,\widetilde{A}x\rangle
$.

It just remains to prove that the drift term $\overline{B}$ is
orthogonal to $K={\rm Ker}\, \overline{A}(y)$. Due to \eqref{hom:c}
and the fact that $K={\rm Ker} \,\overline{A}(y)$ does not depend on
$y\in\R^d$, it suffices to prove that ${\rm Ker}\,
\overline{H}(y)\subset {\rm Ker} \overline{A}(y)=K$. But this is an
easy consequence of \eqref{hom:a}, \eqref{hom:b} and Assumption
\ref{hyp:cont}, especially $|{\bfm H}(\omega,y)|\leq M^2 {\bfm a}
(\omega,y)$.\qed

%%%%%%%%%%%%%%%%%%%%%%%%%%%%%%%%%%%%%%%%%%%%%%%%%%%%%%%%%%%%%%%%%%%%%%%%%%%%%%%%%%%%%%%%%%%%%%%%%%%%

\section{Tightness}\label{sec_tightness}

%%%%%%%%%%%%%%%%%%%%%%%%%%%%%%%%%%%%%%%%%%%%%%%%%%%%%%%%%%%%%%%%%%%%%%%%%%%%%%%%%%%%%%%%%%%%%%%%%%%%%%%%%%%
We now turn to the tightness of the process $X^\varepsilon$, ie we want to prove that the family $(X^\varepsilon)_\varepsilon$ is tight in $C([0,T],\R^d)$ equipped with the uniform topology. That step of our result deeply differs from the uniform
elliptic case \cite{rhodes:06}. Indeed, uniform ellipticity of the
diffusion matrix provides strong transition density estimates of the
process $X^\varepsilon$, the so-called Aronson estimates, from which
the tightness of $X^\varepsilon$ is then easily derived. Of course,
in the degenerate framework, tightness of $X^\varepsilon$ cannot be
tackled this way. The method presented below is inspired from
\cite{wu} and is based on the idea that the process $X^\varepsilon$
is not too far from being reversible at a microscopic scale. The
contributions of the macroscopic variations make a drift appear,
unlike in \cite{wu}.

Let us now go into details. As in Section \ref{sec:auxiliary}, we
can solve the following equation for $i=1,\dots,d$ and $\lambda>0$
\begin{equation}\label{eq:auxtightness}
  \lambda {\bfm w}^i_\lambda(.,y)-{\bfm S}^y{\bfm w}^i_\lambda (.,y)={\bfm b}_i(.,y)
\end{equation}
and get the same properties as in Proposition \ref{correctorlocal},
namely
\begin{proposition}\label{prop:esttight}
For each fixed $y\in \reel^d$ and $1\leq i \leq d $, the family
$(\nabla^{\widetilde{\sigma}} {\bfm w}_\lambda^i(.,y))_\lambda$
converges to a limit $\widetilde{{\bfm \zeta}}_i(.,y)\in
L^2(\Omega)^d$ as $\lambda $ goes to $0$. The same property holds
for the derivatives, that is, the families
$(\nabla^{\widetilde{\sigma}}\partial_{y_j}{\bfm
w}^i_\lambda)_\lambda$, $
(\nabla^{\widetilde{\sigma}}\partial^2_{y_jy_k}{\bfm
w}^i_\lambda)_\lambda$ ($1\leq i,j,k \leq d$) respectively converge
to $\partial_{y_j}\widetilde{{\bfm \zeta}}_i(.,y)$,
$\partial^2_{y_jy_jk}\widetilde{{\bfm \zeta}}_i(.,y)$ in
$L^2(\Omega)^d$. Furthermore, the function ${\bfm w}_\lambda^i$ as
well as its derivatives $\partial_{y_j}{\bfm w}^i_\lambda$,
$\partial^2_{y_jy_k}{\bfm w}^i_\lambda $ satisfy \eqref{eq:convcorr}
and estimates \eqref{eq:borncorr1} and \eqref{eq:borncorr2}, for
some positive constant $C_{\ref{prop:esttight}}$ independent of
$\lambda>0$ and $y\in\reel^d$.
\end{proposition}

As in the proof of Theorem \ref{theoremerg2}, we want to use a time
reversal argument. Once again, we are faced with the lack of
smoothness of $ {\bfm w}_\lambda$ in order to apply the Itô formula.
To overcome this difficulty, we proceed as in Section
\ref{sec:prelim}. Since the arguments are quite similar, we just
outline the main ideas without further details. Let us consider, for
$n\geq 1$, $\lambda>0$ and $1 \leq i \leq d$, the solution $ {\bfm
w}^{i,n}_\lambda$ of the following
 equation
 \begin{equation}\label{eq:auxtightnessn}
  \lambda {\bfm w}^{i,n}_\lambda(.,y)-{\bfm S}^y{\bfm w}^{i,n}_\lambda (.,y)-n^{-1}\Delta{\bfm w}^{i,n}_\lambda (.,y)={\bfm b}_i(.,y)
\end{equation}
Introducing a sequence of regularizing sequence of mollifiers
$(\varrho_m)_{m\in\nat}\in C ^\infty_c(\R^d\times\R^d)$ (smooth
functions with compact support), we define $$ {\bfm w}^{i,n}_
{\lambda,m}(\omega,y)=\int_{\R^{2d}} {\bfm
w}^{i,n}_\lambda(\tau_x'\omega,y-y')\varrho_m(x',y')\,dx'\,dy', $$
which is a smooth function. Following the proof of Theorem
\ref{theoremerg2}, under the invariant measure $e^{-2V(x)}\,dx $ of
the process $X^{n,\varepsilon}$, we can write
\begin{align}\label{eq:tightforw}
w^{i,n}_{\varepsilon^2,m}(\overline{X}^{n,\varepsilon}_t,X^{n,\varepsilon}_t)=&
w^{i,n}_
{\varepsilon^2,m}(\overline{X}^{n,\varepsilon}_s,X^{n,\varepsilon}_s)+
\int_s^t[{\cal L}^{n,\varepsilon} (w^{i,n}_
{\varepsilon^2,m}(\cdot/\varepsilon,\cdot))](\overline{X}^{n,\varepsilon}_r,X^{n,\varepsilon}_r)\,dr\\&+(\overrightarrow{{\cal
M}}^{\varepsilon,n,m}_t-\overrightarrow{{\cal
M}}^{\varepsilon,n,m}_s),\nonumber\\
\label{eq:tightbackw}  w^{i,n}_
{\varepsilon^2,m}(\overline{X}^{n,\varepsilon}_s,X^{n,\varepsilon}_s)=&
w^{i,n}_
{\varepsilon^2,m}(\overline{X}^{n,\varepsilon}_t,X^{n,\varepsilon}_t)+
\int_s^t[({\cal L}^\varepsilon)^*(w^{i,n}_
{\varepsilon^2,m}(\cdot/\varepsilon,\cdot))](\overline{X}^{n,\varepsilon}_r,X^{n,\varepsilon}_r)\,dr\\&+(\overleftarrow{{\cal
M}}^{\varepsilon,n,m}_t-\overleftarrow{{\cal
M}}^{\varepsilon,n,m}_s),\nonumber
\end{align}
where $\overrightarrow{{\cal M}}^{\varepsilon,n,m} $ and
$\overleftarrow{{\cal M}}^{\varepsilon,n,m}  $ are two martingales
respectively with respect to the forward filtration $({\cal
F}^{n,\varepsilon} _s)_{0\leq s\leq T}\equiv
\sigma\left\{X^{n,\varepsilon}_r;0\leq r\leq s\right\}$ and with
respect to the backward filtration $({\cal
G}^{n,\varepsilon}_s)_{0\leq s\leq T}\equiv\sigma\left\{
X^{n,\varepsilon}_r;s\leq r\leq T\right\}$. The quadratic variations
of both martingales match
$$\epsilon^{-2}\int_0^.(D{\bfm w}^{i,n}_
{\varepsilon^2,m}+\varepsilon
\partial_y {\bfm w}
^{i,n}_ {\varepsilon^2,m})^*(a+n^{-1}{\rm Id})(D{\bfm w}^{i,n}_
{\varepsilon^2,m}+\varepsilon
\partial_y {\bfm w}^{i,n}_
{\varepsilon^2,m})(\overline{X}^{n,\varepsilon}_r,X^{n,\varepsilon}_r)\,dr.$$
Adding up \eqref{eq:tightforw} and \eqref{eq:tightbackw}, passing to
the limit as $m\to \infty$ (as explained in \cite[Lemma
5.3]{rhodes:06}) and $n\to \infty$ (as explained in Section
\ref{sec:prelim}) and using (\ref{eq:auxtightness}) leads to
\begin{align}\label{eq:b}
\varepsilon^{-1}\int_s^tb_i(\overline{X}^{\varepsilon}_r,X^{\varepsilon}_r)\,dr
& =\int_s^t[\varepsilon w^i_{\varepsilon^2}+(1/2){\rm
trace}(aD\partial_yw^i_{\varepsilon^2})]
(\overline{X}^{\varepsilon}_r,X^{\varepsilon}_r)\,dr\\
& +\int_s^t\frac{e^{2V}}{2}\big[div_y\big(e^{-2V}a[
Dw^i_{\varepsilon^2}+\varepsilon\partial_y
w^i_{\varepsilon^2}]\big)\big](\overline{X}^{\varepsilon}_r,X^{\varepsilon}_r)\,dr\nonumber\\
& +(1/2)\int_s^t\mathrm{Div}(a)\cdot\partial_yw^i
_{\varepsilon^2}(\overline{X}^{\varepsilon}_r,X^{\varepsilon}_r)\,dr\nonumber\\
&+\varepsilon(\overrightarrow{{\cal M}}^{\varepsilon}_t-\overrightarrow{{\cal
M}}^{\varepsilon}_s)+\varepsilon(\overleftarrow{{\cal
M}}^{\varepsilon}_t-\overleftarrow{{\cal
M}}^{\varepsilon}_s)\nonumber\\& \equiv
E^{1,\varepsilon}_{s,t}+E^{2,\varepsilon}_{s,t}+T^{1,\varepsilon}_{s,t}+T^{2,\varepsilon}_{s,t}\nonumber,
\end{align}
where $\varepsilon\overrightarrow{{\cal M}}^{\varepsilon} $ and
$\varepsilon\overleftarrow{{\cal M}}^{\varepsilon}  $ are two martingales,
respectively with respect to the forward filtration $({\cal
F}^{\varepsilon} _s)_{0\leq s\leq T}\equiv
\sigma\left\{X^{\varepsilon}_r;0\leq r\leq s\right\}$ and with
respect to the backward filtration $({\cal
G}^{\varepsilon}_s)_{0\leq s\leq T}\equiv\sigma\left\{
X^{\varepsilon}_r;s\leq r\leq T\right\}$, with quadratic variations
\begin{equation}\label{eq:quacov}
\int_0^.(D{\bfm w}^{i}_ {\varepsilon^2}+\varepsilon
\partial_y {\bfm w}
^{i}_ {\varepsilon^2})^*a(D{\bfm w}^{i}_ {\varepsilon^2}+\varepsilon
\partial_y {\bfm w}^{i}_
{\varepsilon^2})(\overline{X}^{\varepsilon}_r,X^{\varepsilon}_r)\,dr.\end{equation}
Theorem \ref{theoremerg1} establishes the following  convergence $$
\lim_{\varepsilon\to 0}\em\big[\sup_{0\leq  t\leq
T}\Big|E^{1,\varepsilon}_{0,t}+E^{2,\varepsilon}_{0,t}-\int_0^t\bar{G}(X^\varepsilon_r)\,dr \Big|\big]=0,$$ 
where $$\bar{G}(y)=\M\big[(1/2){\rm trace}({\bfm \sigma}\partial_y{\bfm \xi}_i)(.,y)+(e^{2V}/2) {\rm div}_y\big(e^{-2V}{\bfm \sigma}{\bfm \xi}_i\big)(.,y)\big)\big].$$
From Proposition \ref{prop:esttight} and \eqref{eq:borncorr1}, $\bar{G}$ is bounded so that the tightness
of the process $t\mapsto \int_0^t\bar{G}(X^\varepsilon_r)\,dr$ in $C([0,T],\R) $ results from the Kolmogorov criterion. The tightness of $E^{1,\varepsilon}+E^{2,\varepsilon}$ follows.

Let us investigate now the term $T^{1,\varepsilon}_{s,t}=(1/2)\int_s^t\mathrm{Div}(a)\cdot\partial_yw^i
_{\varepsilon^2}(\overline{X}^{\varepsilon}_r,X^{\varepsilon}_r)\,dr$. Note that it
can not be treated with Theorem \ref{theoremerg1} because the
$L^2$-norm of ${\rm Div}({\bfm a})\partial_y{\bfm
w}_{\varepsilon^2}$ need not be bounded. Inspired by the proof of
Theorem \ref{maintheorem} in Section \ref{sec:proof}, we define
$${\bfm \Psi}^i_\varepsilon\equiv  \mathrm{Div}({\bfm
a})\cdot\partial_y{\bfm w}^i_{\varepsilon^2}+ {\rm trace}({\bfm
a}D\partial_y{\bfm w}^i_{\varepsilon^2})+\varepsilon div_y({\bfm
a})\cdot\partial_y{\bfm w}^i_{\varepsilon^2}+\varepsilon {\rm
trace}({\bfm a}\partial^2_{yy}{\bfm w}^i_{\varepsilon^2})-2\varepsilon{\bfm a}_{pj}\partial_{y_j}{\bfm u}^i_{\varepsilon^2}\partial_{y_p}V.$$ By
making two successive integrations by parts as in Section
\ref{sec:proof}, we establish for any ${\bfm \varphi}\in {\cal
C}\times C^\infty_0(\reel^d)$:
\begin{align*}
\M_V[{\bfm \Psi}^i_\varepsilon,{\bfm \varphi} ] =  -\M_V[{\bfm
a}\partial_y{\bfm w}^i_{\varepsilon^2}\cdot (D{\bfm
\varphi}+\varepsilon\partial_y{\bfm \varphi})]
\stackrel{\mathrm{Prop.} \,\ref{prop:esttight}}{\leq} &
C_\varepsilon\M_V[|{\bfm \sigma}^*(D{\bfm
\varphi}+\varepsilon\partial_y{\bfm \varphi})|^2]^{1/2},
\end{align*}
where the family $(\varepsilon C_\varepsilon)_\varepsilon$ converges
to $0$ as $\varepsilon$ goes to $0$. Theorem \ref{theoremerg2} then
ensures that
$$\em\Big[\sup_{0\leq s \leq
t}\Big(\int_s^t\Psi_\varepsilon(\overline{X}^{\varepsilon}_r,X^{\varepsilon}_r)\,dr\Big)^2\Big]\rightarrow
0$$ as $\varepsilon$ goes to $0$. Thanks to Theorem
\ref{theoremerg1} and Proposition \ref{prop:esttight}, we have
$$\em\Big[\sup_{0\leq s \leq
t}\big|\int_0^s{\rm trace}({\bfm a}D\partial_y{\bfm w}
_{\varepsilon^2})(\overline{X}^{\varepsilon}_r,X^{\varepsilon}_r)\,dr-\int_0^s\bar{\Phi}(X^{\varepsilon}_r)\,dr\big|^2\Big]\to
0$$ as $\varepsilon$ goes to $0$, where
$\bar{\Phi}(y)=\lim_{\varepsilon\to 0}\M[{\rm trace}({\bfm
a}D\partial_y{\bfm w} _{\varepsilon^2})(\cdot,y)] $. The Kolmogorov
criterion and Proposition \ref{prop:esttight} ensure the tightness
 in $C([0,t];\reel)$ of the process $\int_0^\cdot\bar{\Phi}(X^{\varepsilon}_r)\,dr$. Moreover, from Proposition \ref{prop:esttight} and \eqref{eq:inv}, the
 process $\int_0^. \big[\varepsilon div_y({\bfm
a})\cdot\partial_y{\bfm w}^i_{\varepsilon^2}+\varepsilon {\rm
trace}({\bfm a}\partial^2_{yy}{\bfm
w}^i_{\varepsilon^2})-2\varepsilon{\bfm a}_{pj}\partial_{y_j}{\bfm u}^i_{\varepsilon^2}\partial_{y_p}V\big](\overline{X}^{\varepsilon}_r,X^{\varepsilon}_r)\,dr$
converges in law in $ C([0,T];\R)$ to $0$. This proves the tightness
of $T^{1,\varepsilon}$ in $C([0,t];\reel^d)$.

 It just remains to
treat the martingale term $T^{2,\varepsilon}$. According to Theorem
4.13 in \cite{jacod}, it suffices to establish the tightness of the
brackets of these two martingales (see \eqref{eq:quacov}). Their
tightness results from Theorem \ref{theoremerg1}, Proposition
\ref{prop:esttight} and the Kolmogorov criterion again. The
tightness of $X^\varepsilon$ is now clear.  \qed
%%%%%%%%%%%%%%%%%%%%%%%%%%%%%%%%%%%%%%%%%%%%%%%%%%%%%%%%%%%%%%%%%%%%%%%%%%%%%%%%%%%%%%%%%%%%%%%%%%%%%%%%%%%%%%%%%

\end{document}